\documentclass[english,12pt]{smfart}
\usepackage{amsmath,amssymb,amscd,amsfonts,amsmath,mathrsfs,verbatim}
\usepackage[all]{xy}
\usepackage[mathcal]{euscript}
\usepackage{epsfig}
\usepackage{times}
\usepackage{hyperref}

\newcommand{\Z}{{\mathbf Z}}                   
\newcommand{\R}{{\mathbf R}}                   
\renewcommand{\H}{{\mathbf H}}                   
\newcommand{\C}{{\mathbf C}}                   
\newcommand{\A}{{\mathbf A}^1}                   
\renewcommand{\P}{\mathbf{P}^1}                
\newcommand{\PN}{\mathbf{P}^N}               
\newcommand{\E}{{\mathcal E}}                  
\newcommand{\N}{{\mathcal N}}                  
\renewcommand{\O}{{\mathcal O}}                 
\newcommand{\Om}{{\mathbf \Omega}}                 
\newcommand{\Mod}{{\mathcal M}}               
\renewcommand{\d}{\mbox{d}}                      

\newcommand{\isom}{\cong}
\renewcommand{\to}{\longrightarrow}
\newcommand{\dbar}{\bar{\partial}}
\newcommand{\End}{{\mathcal E}nd}         
\renewcommand{\Om}{{\mathbf \Omega}}        
\newcommand{\Hecke}{{\mathcal H}}         
\renewcommand{\L}{{\mathcal L}}           

\DeclareMathOperator{\coker}{coker}       
\DeclareMathOperator{\im}{Im}             
\DeclareMathOperator{\tr}{tr}             
\DeclareMathOperator{\Ext}{Ext}           
\DeclareMathOperator{\diag}{diag} 
\DeclareMathOperator{\res}{res} 

\renewcommand{\theenumi}{(\roman{enumi})}

\newtheorem{prop}{Proposition}[section]
\newtheorem{cond}[prop]{Condition}
\newtheorem{rk}[prop]{Remark}

\newtheorem{lem}[prop]{Lemma}

\newtheorem{cor}[prop]{Corollary}
\newtheorem{thm}[prop]{Theorem}
\newtheorem{qn}[prop]{Question}

\title{Deformations of Fuchsian equations and logarithmic connections}
\author[Szil\'ard Szab\'o, R\'enyi Institute of
  Mathematics, Budapest]{Szil\'ard Szab\'o}
\date{\today}

\begin{document}

\begin{abstract}
We give a geometric proof to the classical fact that the dimension of the deformations of a given 
generic Fuchsian equation without changing the semi-simple conjugacy class of its local monodromies 
(``number of accessory parameters'') is 
equal to half the dimension of the moduli space of deformations of the associated local system. 
We do this by constructing a weight $1$ Hodge structure on the infinitesimal deformations of 
logarithmic connections, such that deformations as an equation correspond to the $(1,0)$-part. 
This answers a question of Nicholas Katz, who noticed the dimension doubling mentioned above. 
We then show that the Hitchin map restricted to deformations of the Fuchsian equation is a 
one-to-one \'etale map. 
Finally, we give a positive answer to a conjecture of Ohtsuki about the maximal number of apparent 
singularities for a Fuchsian equation with given semisimple monodromy, and define a Lagrangian 
foliation of the moduli space of connections whose leaves consist of logarithmic connections 
that can be realised as Fuchsian equations having apparent singularities in a prescribed finite set. 
\end{abstract}

\maketitle

\section{Introduction}
Let $p_1,\ldots ,p_n\in \P$ be $n \geq 2$ fixed points in the affine part of the complex projective line, 
and let $p_{0}$ be the point at infinity. 
Define $P$ to be the simple effective divisor $p_0 + \cdots + p_n$ in $\P$, and let $P^o=p_1 + \cdots + p_n$. 
Consider the function 
\begin{equation}\label{psi}
    \psi(z) = (z-p_1)\cdots (z-p_n)
\end{equation}
as an identification between $\O_{\P}$ and $\O_{\P}(-P^o)$ on the affine $\A=\P\setminus \{p_0 \}$. 

Let $w=w(z)$ be a holomorphic function of the complex variable $z$, $w^{(k)}$ its $k$-th order differential 
with respect to $z$ and 
\begin{equation}\label{equation}
   w^{(m)}-\frac{G_1(z)}{\psi}w^{(m-1)}-\cdots -\frac{G_m(z)}{\psi^m}w=0,
\end{equation}
where the $G_k$ are polynomials in $z$, be a Fuchsian differential equation. 
We recall that this means that all the solutions $w$ grow at most polynomially 
with $(z-p_j)^{-1}$ near $p_j$ for any $1 \leq j \leq n$ (respectively, with $z$ near infinity). 
The left-hand side of (\ref{equation}) is a linear differential operator of the order $m$ of $w$, 
that we shall denote by $\L$. We shall identify the equation (\ref{equation}) with the operator $\L$. 
Fuchs gave a necessary and sufficient condition on the degrees of the $G_k$'s for $\L$ to 
be of Fuchsian type: the degree of $G_k$ has to be at most $k(n-1)$. 

Let us introduce the expressions 
\begin{align}
    w_1 & = w \notag \\
    w_2 & = \psi \frac{\d w}{\d z} \notag \\
        & \vdots \label{extension} \\
    w_m & = \psi^{m-1} \frac{\d^{m-1} w}{\d z^{m-1}}\notag . 
\end{align}
We think of $w$ as a local holomorphic function in $z$, or in modern terminology a local section 
of the structure sheaf $\O$. Then, on the affine open $\A$ the function $w_k$ is meromorphic 
with zeroes of order at least $(k-1)$ in $P$. 
In other words, on the affine part $\A$ a vector $(w_1,\ldots,w_m)$ is a section of the 
holomorphic vector bundle 
$$
    \tilde{E}= \O_{\P} \oplus \O_{\P}(-P^o) \oplus \ldots \oplus \O_{\P}((1-m)P^o). 
$$ 
We equip $\tilde{E}$ with the algebraic integrable connection with logarithmic poles at $P$ 
\begin{equation}\label{conn}
    D_{\L}=\d^{1,0}-\frac{A(z)}{\psi(z)} \d z, 
\end{equation}
where $A(z)$ is the {\em modified companion matrix} 
\begin{equation}\label{a}
    A  =\begin{pmatrix}
             0 & 0& 0 & 0 &\ldots \ldots & G_m \\
          1 & \psi ' & 0& 0 & \ldots \ldots &G_{m-1} \\
          0 & 1 & 2 \psi ' & 0 & \ldots \ldots &G_{m-2} \\
          \vdots  & \vdots & \vdots &   &  & \vdots  \\
          0 & 0  & 0 & \ldots   &  (m-2)\psi '& G_2 \\
          0 & 0 & 0 & \ldots & 1 & G_1 + (m-1)\psi ' 
      \end{pmatrix}
\end{equation}
of equation (\ref{equation}). Here we have denoted $\psi '=\d \psi /\d z$. 
One readily checks that a meromorphic function $w$ on some open set $U\subset \A$ with poles at most in $P$ 
locally solves (\ref{equation}) if and only if the vector $(w_1=w,w_2,\ldots ,w_m)$ is a parallel section of 
$D_{\L}$ on $\tilde{E}$ over $U$ for some (hence, only one) vector $(w_2,\ldots ,w_m)$. 

\begin{rk}
Instead of the formulae (\ref{extension}), we could have simply chosen
$w_k=w^{(k)}$, and the form of the same connection $D_{\L}$ in this
trivialisation would then be a usual companion matrix. The reason for 
our choice for the extension (\ref{extension}) is that it gives rise to 
a logarithmic lattice. In fact, the two points of view are equivalent, so 
that a connection in modified companion form is also induced by an equation. 
\end{rk}

By assumption, the integrable connection $D_{\L}$ is regular at infinity as well. 
Therefore, according to a theorem of 
N. Katz (Thm. II.1.12 \cite{del}), there exists a lattice for the meromorphic bundle 
$$
  \N = \O_{\P}(*P) \oplus \ldots \oplus \O_{\P}(*P)
$$ 
at infinity with respect to which $D_{\L}$ is a logarithmic connection, 
i.e. its local form in any holomorphic trivialisation contains $1$-forms with at most first-order poles. 
Such a lattice at infinity can be obtained similarly to the case of the other singularities. 
For this purpose, let $\zeta=z^{-1}$ be a local coordinate at infinity. Recall that the 
first component $w=w_1$ of (\ref{extension}) is supposed to be a section of $\O_{\P}$. Then, a logarithmic 
lattice $(\tilde{w}_1,\ldots,\tilde{w}_m)$ at infinity can be obtained by the formulae 
\begin{align}
    \tilde{w}_1 & = w \notag \\
    \tilde{w}_2 & = \zeta \frac{\d w}{\d \zeta} \notag \\
        & \vdots \label{extinf} \\
    \tilde{w}_m & = \zeta^{m-1} \frac{\d^{m-1} w}{\d \zeta^{m-1}}, \notag 
\end{align}
and the form of the connection is then again a modified companion matrix with entries in the 
last column equal to the coefficients of the equation multiplied by an appropriate power of $\zeta$. 
Now, as 
$$
    \zeta^{j} \frac{\d^{j} w}{\d \zeta^{j}} = (-z)^{j} \frac{\d^{j} w}{\d z^{j}} 
$$
and for large $z$ one has 
$$
    \psi(z) \approx z^n, 
$$
we deduce that the trivialisations (\ref{extension}) and (\ref{extinf}) are linked on $\C^*$ by the 
matrix 
$$
    \diag(1,-z^{n-1},\ldots,(-1)^{m-1}z^{(m-1)(n-1)}). 
$$
It follows that the connection $D_{\L}$ extends logarithmically onto the holomorphic bundle 
\begin{equation}\label{holbdl}
    E_{\L} = \O_{\P} \oplus \O_{\P}(\infty - P^o) \oplus \ldots \oplus \O_{\P}((m-1)(\infty - P^o))
\end{equation}
over $\P$. 
Hence, we get a bijective correspondence between local solutions of (\ref{equation}) and 
local parallel sections of the logarithmic connection $D_{\L}$ on the holomorphic bundle $E_{\L}$. 
The following fundamental result is part of folklore; however, we give a proof here for lack 
of appropriate reference. 
\begin{prop}\label{prop:gaugeequiv}
If the logarithmic connections induced by the Fuchsian equations $\L_1,\L_2$ are gauge-equivalent, 
then $\L_1=\L_2$. 
\end{prop}
\begin{proof}
Suppose there exists a gauge transformation $g \in \End(E_{\L_1},E_{\L_2})$ mapping $D_{\L_1}$ into $D_{\L_2}$. 
In the trivialisations (\ref{holbdl}) of $E_{\L_1}$ and $E_{\L_2}$, $g$ can be written as a matrix whose 
entry in the $k$-th row and $l$-th column is a global holomorphic section of the sheaf $\O((k-l)(n-1))$. 
It follows that the matrix of $g$ is lower triangular, and that the entries on the diagonal are 
global sections of the trivial holomorphic line bundle over $\P$, hence constants. 
For $j=1,2$ let us write on the affine part $\C$ of $\P$ away from infinity the expressions 
$$
    D_{\L_j}=\d^{1,0}-\frac{A_j(z)}{\psi(z)} \d z.  
$$
It is then a well-known fact that the action of $g$ on $D_{\L_1}$ is 
$$
   g \cdot (D_{\L_1}) = \d^{1,0} - \frac{g^{-1}A_1(z)g}{\psi(z)} \d z 
     - g^{-1}d^{1,0}g.
$$
It follows from the above observations that $g^{-1}d^{1,0}g$ is strictly lower triagular.  
In particular, the entries on and above the diagonal in the matrices $A_1$ and $A_2$ must agree. 

Let us first consider the case $m=2$. Here, one has 
$$
   g = \begin{pmatrix}
       g_{11} & 0 \\
       g_{21} & g_{22}
       \end{pmatrix}, 
$$
where $g_{21}$ is a global section of $\O(n-1)$ and $g_{11},g_{22}$ are constants with $g_{11}g_{22}\neq 0$, 
and the inverse of this matrix is 
$$
   g^{-1} = \frac{1}{g_{11}g_{22}} 
       \begin{pmatrix}
       g_{22} & 0 \\
       -g_{21} & g_{11}
       \end{pmatrix}. 
$$
Furthermore, the matrices of the equations are 
$$
   A_j= \begin{pmatrix}
       0 & 1 \\
       G^j_2 & G^j_1 + \psi'
     \end{pmatrix}, 
$$
where $G^j_2, G^j_1$ are the coefficients of $\L_j$. 
Straightforward matrix multiplication yields 
$$
   g^{-1}A_1(z)g = \frac{1}{g_{11}g_{22}} 
    \begin{pmatrix}
      g_{21}g_{22} & g_{22}^2 \\
      * & * 
       \end{pmatrix}. 
$$
By the above, the terms in the first row must be equal to $0$ and $1$ respectively. 
We infer that $g_{21}=0$ and $g_{22}=g_{11}$, hence $g$ is a multiple of the identity. 

We now come to the general case. As the computations are more involved but 
of the same spirit, we only sketch the proof. The matrix $g$ is equal to 
$$
   g = \begin{pmatrix}
       g_{11} & 0 & \ldots & 0 \\
       g_{21} & g_{22}  & \ldots & 0 \\ 
       \vdots & & \ddots & \vdots \\
       g_{m1} & g_{m2} & \ldots & g_{mm} 
       \end{pmatrix}, 
$$
and its inverse is of the form 
$$
    g^{-1}=\frac{1}{g_{11}g_{22}\cdots g_{mm}}
      \begin{pmatrix}
      g_{22}\cdots g_{mm} & 0 & \ldots & 0 \\
      * & g_{11}g_{33}\cdots g_{mm}& \ldots & 0 \\
      \vdots & & \ddots & \vdots \\
      * & * & \ldots & g_{11}\cdots g_{m-1,m-1} 
      \end{pmatrix}. 
$$
One has 
\begin{align*}
   g^{-1}A_1(z)g & = \frac{1}{g_{11}g_{22}\cdots g_{mm}} \cdot \\
   & \begin{pmatrix}
      g_{22}\cdots g_{mm} & 0 & \ldots & 0 \\
      * & g_{11}g_{33}\cdots g_{mm}& \ldots & \vdots \\
      \vdots & & \ddots & 0 \\
      * & * & \ldots & g_{11}\cdots g_{m-1,m-1} 
   \end{pmatrix} \\ 
   & \begin{pmatrix}
       g_{21} & g_{22} & 0 & \ldots & 0 \\
       g_{21}\psi'+g_{31} & g_{22}\psi'+g_{32} & g_{33} & \ldots & \vdots \\ 
       \vdots & & \ddots & \ddots & 0 \\
       * & & & g_{m-1,m-1}\psi'+g_{m,m-1} & g_{mm}\\
       * & \ldots & & * & * 
   \end{pmatrix}.
\end{align*}
As the entries on and above the diagonal in this product have to be equal to those of the modified 
companion matrix $A_2(z)$, we deduce as before that $g_{11}=g_{22}=\cdots =g_{mm}$ and 
$g_{21}=g_{32}=\cdots =g_{m,m-1}=0$. 
It follows that right below the diagonal all the entries of the matrix $g^{-1}d^{1,0}g$ vanish. 
Considering now the first sub-diagonal in the product above, we deduce that $g_{31}=g_{42}=\cdots=g_{m,m-2}=0$.
It follows that on the second sub-diagonal of the matrix $g^{-1}d^{1,0}g$ all the entries vanish. 
Continuing this argument, we eventually obtain that all the $g_{kl}$ for $k>l$ must vanish. 
This concludes the proof. 
\end{proof}

Consider now the residue $\res(p,D_{\L})$ of $D_{\L}$ on $E_{\L}$ at 
each of the singular points $p \in P$; it is a well-defined endomorphism of the fiber of $E_{\L}$ at $p$. 
Denote by $\mu^j_1,\ldots,\mu^j_m$ the eigenvalues of the residue at $p_j$. 

\begin{rk}
Notice in particular that $\deg(E_{\L})=(1-n)m(m-1)/2$, in accordance with the classical Fuchs' relation, 
which states that the sum $\sum_{j,k}\mu^j_k$ of the eigenvalues of the residues of $D_{\L}$ in all 
singularities (including infinity) is equal to $(n-1)m(m-1)/2$. 
\end{rk}

Throughout the paper, we will assume the genericity conditions : 
\begin{cond}\label{cond}
The eigenvalues $\mu^j_1,\ldots,\mu^j_m$ of the residue of the integrable connection 
$D_{\L}$ in each singularity $p_j$ do not differ by integers. 
(In particular, they are distinct.) We call this the \emph{non-resonance} condition. 
Furthermore, for any $1 \leq k <m$ there exists no choice of $k$-tuples of eigenvalues 
$\mu^j_{l^j_1},\ldots,\mu^j_{l^j_k}$ at all singular points such that 
$\sum_{j=0}^m\sum_{r=1}^k\mu^j_{l^j_r} \in \Z$. 
\end{cond}

\begin{rk}
The second condition implies that any logarithmic integrable connection $(E,D)$ with these residues 
is stable in the usual sense: any $D$-invariant subbundle of $E$ has slope smaller than $E$. 
Indeed, by the residue theorem for such a logarithmic connection there exist no non-trivial 
$D$-invariant subbundles at all. 
Stability is needed at two instances: for the definition of the moduli space $\Mod$ of integrable 
connections containing $D_{\L}$ (see Remark \ref{rk:moduli}), and to be able to apply non-Abelian Hodge 
theory in Section \ref{sec:Hitchinmap}. 
\end{rk}

We are interested in the following two numbers: 
\begin{enumerate}
\item the dimension $e$ of the space $\E$ of deformations of the polynomials in (\ref{equation}) 
so that all residues of the associated integrable connection $D_{\L}$ remain in the same conjugacy class 
\label{def1}
\item the dimension $c$ of the moduli space $\Mod$ of S-equivalence classes of (semi-)stable integrable 
connections $(E,D)$ logarithmic in $P$ over a vector bundle $E$ of degree $d=(1-n)m(m-1)/2$, 
with fixed conjugacy classes of all its residues. \label{def2} 
\end{enumerate}
\begin{rk}\label{rk:moduli}
For the definition of the moduli space $\Mod$, see for example Sections 6-8 of \cite{biqboa}. 
An alternative way would be to define it geometrically as a symplectic leaf of the coarse moduli 
scheme of stable logarithmic connections (with arbitrary residues) constructed in Theorem 3.5 of \cite{nit}. 
On the other hand, the space $\E$ is well-known to be an affine space, see e.g. \cite{ince}. 
\end{rk}
It is immediate that $e\leq m$, for the space of the Fuchsian deformations \ref{def1} is contained in the 
space of integrable connections \ref{def2} having the right monodromy, and two integrable connections induced 
by different Fuchsian equations cannot be gauge-equivalent, so this inclusion map is injective. 
In short, we will call deformations leaving invariant the conjugacy classes of all the residues 
\emph{locally isomonodromic}.
Using Fuchs' condition, the number $e$ was computed by Forsyth in \cite{for}, pp. 127-128. 
In the introduction of his book \cite{ka}, N. Katz computed $c$, and noticed that 
\begin{equation}\label{equality}
    c=2e. 
\end{equation}
In fact, both sides of this equation turn out to be 
\begin{equation}\label{exactvalue}
   2 - 2m^2 + {m(m-1)}(n+1). 
\end{equation}
He also asked whether a geometric reason underlies this equality, more precisely, whether a weight $1$ 
Hodge structure can be found on the tangent to the moduli space of integrable connections, whose $(1,0)$-part  
would give precisely the locally isomonodromic deformations of Fuchsian equations. 

We will define such a Hodge structure in Section \ref{sec:Hodgestr}: 
\begin{thm}\label{thm:main} 
Let $D_{\L}$ be an integrable connection (\ref{conn})-(\ref{a}) induced by a Fuchsian equation (\ref{equation}) 
satisfying Condition \ref{cond}. Then there exists a natural weight $1$ Hodge structure on the tangent 
at $D_{\L}$ to the moduli space $\Mod$ of integrable connections
$$
     T_{D_{\L}}\Mod = H^{1,0} \oplus H^{0,1} 
$$
with the property that its part of type $(1,0)$ is the tangent of the space $\E$ of locally isomonodromic 
deformations of the Fuchsian equation: 
$$
    T_{D_{\L}}\E = H^{1,0}.
$$ 
\end{thm}
The Hodge structure comes from a hypercohomology long exact sequence, and is well-defined on 
the tangent space at all elements of the moduli space. A similar exact sequence already 
appears in Proposition 4.1 of \cite{nit}. 

In Sections \ref{sec:Hitchinmap} and \ref{sec:Applications}, we show: 
\begin{thm}\label{thm:Hitchinmap}
The Hitchin map restricted to $\E$ is a bijection. 
Furthermore, $\E$ is an algebraic subvariety which is Lagrangian with respect to the 
natural holomorphic symplectic structure of the de Rham moduli space. 
\end{thm}
For precise definitions, see Section \ref{sec:Hitchinmap}. The Lagrangian property is also 
proved independently by J. Aidan \cite{aidan}. 

Finally, in Corollary \ref{cor:Lagrfol} we use a result of Ohtsuki to determine the exact number 
of apparent singularities of the Fuchsian equation associated to a generic logarithmic connection, 
and deduce that $\Mod$ is foliated in Lagrangian subspaces: 
\begin{thm}
The smallest number $N$ such that any logarithmic connection $(E,D) \in \Mod$ corresponds 
to a Fuchsian equation with at most $N$ apparent singularities is equal to $e$. 
The subspaces of $\Mod$ consisting of connections that can be realised by a Fuchsian 
equation with fixed locus of apparent singularities define a Lagrangian foliation of $\Mod$. 
\end{thm}
The first statement of this theorem is proved independently by B. Dubrovin and M. Mazzocco in \cite{dubma}.

\section{The Hodge structure} \label{sec:Hodgestr}

\subsection{Construction}\label{subsec:construction} 
Let $(E,D)$ be an arbitrary element of $\Mod$. 
Let us denote by $\End_{iso}(E)$ the sheaf of locally isomonodromic endomorphisms, 
which are by definition the  endomorphisms $\varphi$ whose value $\varphi(p_j)$ at $p_j$ 
lies in the adjoint orbit of the residue of the connection in the Lie algebra $\mathfrak{gl}(m)$. 
By Condition \ref{cond}, this residue is regular diagonal in a suitable trivialisation; 
the locally isomonodromic endomorphisms are the ones whose value at $p_j$ is an off-diagonal matrix 
in this basis. 
The infinitesimal deformations of the integrable connection $D$ (without changing the eigenvalues of the 
residues) are then described by the first hypercohomology $\H^1(D)$ of the two-term complex 
\begin{equation}\label{complex}
    \End(E) \xrightarrow{D} \Om^1(P) \otimes \End_{iso}(E)
\end{equation}
(see Section 12 of \cite{biq}). 

Denoting by $\H^i$ the $i$-th hypercohomology of a complex, the hypercohomology long exact sequence 
for (\ref{complex}) reads 
\begin{align}\label{exseq}
   0 \to \H^0(D) & \to H^0(\End(E)) \xrightarrow{H^0(D)} H^0(\Om^1(P) \otimes \End_{iso}(E)) \to \notag \\
     \to \H^1(D) & \to H^1(\End(E)) \xrightarrow{H^1(D)} H^1(\Om^1(P) \otimes \End_{iso}(E)) \to \notag \\
     & \to \H^2(D) \to 0. 
\end{align}
The maps $H^i(D)$ are induced by $D$ on the corresponding cohomology spaces. 
Setting 
\begin{align}
    C & = \coker(H^0(D)) \\
    K & = \ker(H^1(D)), 
\end{align}
there follows a short exact sequence for the space of infinitesimal deformations: 
$$
   0 \to C \to \H^1(D) \to K \to 0. 
$$
The term $C$ roughly corresponds to infinitesimal modifications of the $(1,0)$-part of the integrable 
connection while keeping the holomorphic structure fixed, 
whereas the term $K$ corresponds to infinitesimal modifications of the holomorphic structure. 

\begin{lem}\label{lem:duality}
We have $K^{\vee}\cong C$.
\end{lem}
\begin{proof}
Let us first compute the dual of complex (\ref{complex}). 
The dual of the sheaf $\End(E)$ is clearly $\End(E)$ itself, so the element of degree $1$ 
in the dual is $\Om^1 \otimes\End(E)$. 
As $\End_{iso}(E)$ is the sheaf of endomorphisms of $E$ vanishing on the diagonal in a 
local diagonalising trivialisation of $D$ near the singularities, and the dual of 
the vanishing condition is having a simple pole, it follows that the element $\Hecke \End(E)$ 
of degree $0$ in the dual complex fits into the short exact sequence 
$$
   0 \to \Hecke \End(E) \to \End(E) \to \bigoplus_{p \in P} \im_p \to 0, 
$$
where $\im_p$ stands for the skyscraper sheaf supported at $p$, with stalk equal to 
the off-diagonal part of $\End(E)_p$ in a diagonalising trivialisation, and the map to 
$\im_p$ is evaluation followed by projection to the off-diagonal part. 
This stalk is just the image at $p$ of the adjoint action of the residue 
of the connection on endomorphisms, so it is intrinsically defined. 
We will call $\Hecke \End(E)$ the Hecke transform of $\End(E)$ along the image of the residue. 
One also checks immediately that the terms of degree $1$ of (\ref{complex}) and its dual 
fit to a similar short exact sequence: 
$$
   0 \to \Om^1 \otimes \End(E) \to \Om^1(P) \otimes \End_{iso}(E) \to 
      \bigoplus_{p \in P}\im_p \to 0, 
$$
where the map to $\im_p$ is taking residue at $p$ followed by projection to the image of the residue 
of $D$. In other words, the dual of the complex (\ref{complex}) is linked to (\ref{complex}) 
by a Hecke transformation along the image of the residue. The hypercohomology long exact sequence 
(\ref{exseq}) is the long exact sequence associated to a short exact sequence of two-term 
complexes. All non-zero terms of these two-term complexes for (\ref{complex}) and its dual 
are related by Hecke transformation along the image of the residue. As the residue of $D$ 
acts by an isomorphism on its image, it follows that the connecting morphism of the long 
exact sequence of the skyscraper sheaves is an isomorphism. As the hypercohomology long exact sequence 
is functorial, we deduce that the spaces 
$$
   \coker(H^0(\Hecke \End(E))\xrightarrow{H^0(D)}H^0(\Om^1 \otimes \End(E))) 
$$
and 
$$
   \coker(H^0(\End(E))\xrightarrow{H^0(D)}H^0(\Om^1(P) \otimes \End_{iso}(E))) 
$$
are isomorphic. We conclude by Verdier duality. 
\end{proof}
Hence, Lemma \ref{lem:duality} together with (\ref{exseq}) exhibits $\H^1(D)$ 
as the extension of two vector spaces of the same dimension: 
$$
    0 \to C \to \H^1(D) \to K \to 0. 
$$ 
Let us show that the two factors of this extension are complex conjugate to each other. 
By equality of the dimensions, it is sufficient to construct an injective anti-linear map 
from one space to the other. 
An element of $C$ is a class of global sections of $\Om^1 (P) \otimes \End_{iso}(E)$, 
modulo the image of the map $H^0(D)$ induced by the connection on global sections of $\End(E)$. 
Let $\alpha' \d z$ represent such a class. 
Denote by $h$ the harmonic metric for the integrable connection and by $^*$ the operation of taking 
adjoint with respect to $h$. 
Then, one has $\coker(H^0(D))\isom \ker(H^0(D)^*)$, therefore there exists a unique global section 
$\alpha\d z \in \ker(H^0(D)^*)$ such that $\alpha\d z=\alpha'\d z + D f$ for some 
$f \in H^0(\End(E))$. 
Write $\nabla=D+\dbar^E$ for the differential geometric flat connection associated to the 
couple $(E,D)$. We denote by $\d_{\nabla}$ the exterior derivative induced by $\nabla$ and by 
$\Delta=\d_{\nabla}\d_{\nabla}^*+\d_{\nabla}^*\d_{\nabla}$ the associated Laplace operator. 
Then, since $\alpha \d z$ is of type $(1,0)$, it follows that 
\begin{align*}
   \Delta(\alpha\d z) & = \d_{\nabla}^* (\dbar^E \alpha\d z)+ \d_{\nabla} (D^*\alpha\d z) \\
       & = \nabla (D^*\alpha\d z) \\
       & = D (D^*\alpha\d z) \\
       & = 0. 
\end{align*}
The second equality in this sequence holds because $\alpha\d z$ is by assumption a 
global holomorphic form, the third one follows from the map 
$$
   H^0(D)^* : H^0(\Om^1 (P) \otimes \End_{iso}(E)) \to H^0(\End(E)), 
$$
and the last one is a consequence of the assumption $\alpha\d z \in \ker(H^0(D)^*)$. 
Hence, $\alpha\d z$ is an endomorphism-valued global harmonic $(1,0)$-form representing the class 
of $\alpha'\d z$. 
Take the adjoint of the endomorphism with respect to $h$ and the complex conjugate $1$-form: 
$\alpha^*\d \bar{z}$. 
According to the $L^2$ Dolbeault resolution (Lemma 9.1 and Theorem 5.1, \cite{biq}) 
it defines an element $[\alpha^*\d \bar{z}]$ of $H^1(\P,\End(E))$. 
Since $\Delta$ is a real operator, we deduce that 
\begin{align*}
   0=\Delta(\alpha^*\d \bar{z}) & =\d_{\nabla}^* (D \alpha^*\d \bar{z}) + 
         \nabla ((\dbar^E)^* \alpha^*\d \bar{z}) \\
      & = \d_{\nabla}^* (D \alpha^*\d \bar{z}), 
\end{align*}
where the second equality follows from $\dbar^E \alpha\d z=0$. 
In particular, one has 
$$
    D \alpha^*\d \bar{z} \in \ker((\dbar^E)^*), 
$$
therefore $D \alpha^*\d \bar{z}$ represents a class in $H^1(\P,\Om^1 (P) \otimes \End_{iso}(E))$. 
Furthermore, $D \alpha^*\d \bar{z}$ is a global holomorphic $2$-form with values in 
$\Om^1 (P) \otimes \End_{iso}(E)$ satisfying $D^*(D \alpha^*\d \bar{z})=0$. 
By irreducibility of $D$, this $2$-form must be trivial, and in particular the cohomology class 
$[\alpha^*\d \bar{z}] \in H^1(\P,\End(E))$ is in the kernel of $H^1(D)$, that is, in $K$. 
Call this cohomology class $\varsigma([\alpha\d z])$. Clearly, $\varsigma$ is then an anti-linear map 
from $C$ to $K$. In the next paragraph we will see that it is bijective. 

The map $\varsigma$ can be extended to define a conjugation on the whole tangent 
space $\H^1(D)$. 
Indeed, by Theorem 12.6 of \cite{biq} the hypercohomology $\H^1(D)$ is isomorphic to 
the $L^2$-cohomology of $\nabla$ with respect to $h$, or said differently, to the 
$L^2$-kernel of the Laplacian $\Delta$. 
Let $\alpha\d z +\beta\d \bar{z}$ be the endomorphism-valued $L^2$ harmonic $1$-form 
representing a given class $a\in\H^1(D)$. 
Since the Laplacian is a real operator, the $1$-form $\beta^*\d z +\alpha^*\d \bar{z}$ 
is also harmonic, hence it represents a class in $\H^1(D)$. We declare this class 
to be $\varsigma(a)$. This map is clearly involutive, in particular bijective. 
The vectors fixed by $\varsigma$ are precisely the harmonic $1$-forms of the form 
$\alpha\d z +\alpha^*\d \bar{z}$. Call $\H^1_{\R}(D)$ the vector space of such forms. 
Then, we have $\H^1(D)=\H^1_{\R}(D)\otimes_{\R} \C$. 

The conjugation $\varsigma$ induces a Hermitian metric $g$ on $\H^1(D)$. 
Let $\alpha_j\d z +\beta_j\d \bar{z}$ be the harmonic representatives of the classes $a_j\in \H^1(D)$ 
for $j \in \{1,2\}$, then it is defined by the formula 
\begin{align}\label{Hermitian}
    g(a_1,a_2)& = a_1 \cup \varsigma(a_2) \\   
  & = \int_{\P}\tr(\alpha_1\d z\wedge \alpha_2^*\d \bar{z}+\beta_1\d \bar{z}\wedge\beta_2^*\d z). \notag
\end{align}
Clearly, $C$ and $K$ are orthogonal complements of each other with respect to this metric. 
In particular, we deduce the orthogonal decomposition of vector spaces 
\begin{align}\label{decomp}
     \H^1(D) & = C \oplus K \\
              & = H^{1,0} \oplus H^{0,1}. \notag
\end{align}
which is just the decomposition of harmonic $1$-forms according to type. 
The second line of this decomposition together with the real structure $\H^1_{\R}(D)$ 
define a weight $1$ Hodge structure on the tangent of $\Mod$, and the Hermitian 
metric (\ref{Hermitian}) induces a polarisation on it. 
In the next subsection we show that it admits the desired property. 

\subsection{Characterisation of deformations of the Fuchsian equation}\label{subsec:char}
In this subsection, we show that the Hodge structure defined in the previous subsection satisfies 
the property claimed in Theorem \ref{thm:main}. In all this part, we consider a point 
$(E_{\L},D_{\L})\in \E$, induced by the Fuchsian equation $\L$. 

We start by showing that infinitesimal deformations of the Fuchsian equation lie inside the $(1,0)$-part. 
Let $\L(t)$ be a deformation of the Fuchsian equation (\ref{equation}) for $t$ in a small disk around 
$0$ in $\C$; by this we mean that we are given polynomials $G_1(t),\ldots G_m(t)$ in the variable $z$ 
with $G_k(0)=G_k$, such that $\deg(G_k(t))\leq k(n-1)$ for each $t$. 
It is clear that the lattices (\ref{extension}) and (\ref{extinf}) are independent of $t$, just as 
the gluing matrices between them. 
It follows that the underlying vector bundle $E_{\L(t)}$ of the connection associated to $\L(t)$ is 
constant. Since the map 
$$
    \H^1(D) \to H^1(\End(E)) 
$$
in the hypercohomology long exact sequence is restriction of an infinitesimal deformation of 
logarithmic connections to the infinitesimal 
deformation of the underlying holomorphic bundle, we deduce that infinitesimal deformations of the 
Fuchsian equation form a subspace of $C$. 

We now come to surjectivity. 
As the spaces $H^{1,0}(D)$ clearly change continuosly with the connection $D$, any vector of 
$H^{1,0}(D_{\L})$ can be extended to a vector field $V$ tangent to $H^{1,0}$ in a small neighborhood of 
$(E_{\L},D_{\L})$ in $\Mod$. 
This vector field can then be integrated to define a flow $(E(t),D(t))$ for $|t|<\varepsilon$ 
in $\R$, with $D(0)=D_{\L}$. 
We will then show that for all such deformations, we have $(E(t),D(t))\in \E$ for all $t$. 
By a result of N. Katz \cite{ka2}, for each $t$ there exists Zariski locally on $\P$ a cyclic vector 
for $(E(t),D(t))$, that is to say away from some finite set $S=S(t)$ of $\P$ a vector $v=v(t)$ in $E$ such 
that $v,D(t)v,\ldots ,D(t)^{m-1}v$ generate $E$. 
Put in different terms, this means that in this cyclic trivialisation of $E(t)$ on 
$\P \setminus (P \cup S)$ the matrix of $D(t)$ is in companion matrix-form. 
Still in other terms, $D(t)$ is the integrable connection 
associated to the Fuchsian differential equation whose coefficients are the functions in 
the lowest row of this companion matrix. The only problem is that the set $S$ might be strictly larger 
than the singular set $P$. The additional points $S\setminus P$ are usually called 
{\em apparent singularities}. In these points, all solutions of the local system associated to $D(t)$ 
extend smoothly, although in the cyclic trivialisation the local form of the integrable connection has 
singularities. 

An element in $H^{0,1}(D)$ can be given by the data of sections 
$h_{ij} \in \Gamma(U_i\cap U_j,\End_{par}(E))$ for some cover $(U_i)_i$ of $\P$ by finitely many open sets,
satisfying the cochain-relation. We choose the open cover in such a way that every singular point 
is in the interior of exactly one open set. 
Such an element is equal to $0$ in cohomology if and only if  there exist sections 
$h_i\in \Gamma(U_i,\End_{par}(E))$ such that $h_i-h_j=h_{ij}$. Suppose we have such a deformation 
at any $t$. Then, possibly passing to a smaller $\varepsilon>0$ and 
slightly restricting the open sets $U_i$, there exists a continuous family in $t$ of sections 
$h_i(t)\in \Gamma(U_i,\End_{par}(E))$ with the property that $h_i(0)=\mbox{Id}$ and that $h_i(t)$ links the 
cyclic vector at time $t$ to the cyclic vector at time $0$ in $U_i$. 
Our aim is to show that no apparent singularities can occur for $D(t)$ changing in this family. 

In order to prove this statement, first we need to understand where the apparent singularities might 
possibly appear. 
The construction of a cyclic vector in \cite{ka2} is clearly holomorphic with respect to the connection 
$D$. Therefore, the locus where it admits a singularity depends holomorphically on the connection $D$ 
in the whole moduli space $\Mod$. 
In particular, this is the case with the positions of the apparent singularities. 
Therefore, the positions of the apparent singularities form the branches of a multi-valued meromorphic 
function $\tau$ on $\Mod$. 
By assumption, there are no apparent singularities for $D(0)=D_{\L}$. This implies that $\tau$ at $D(0)$ 
takes its values in the set $P$. Hence, if some apparent singularities appear for $D(t)$, 
they can only split off from some of the real singularities $p\in P$. 
Suppose $a(t)$ is such a branch splitting off $p$, i.e. such that $a(0)=p$, but not identical to $p$. 
For simplicity of the explanations, assume moreover that $a(t)$ is the only branch splitting off $p$. 
(The case of a branch splitting off $\infty$ is identical, for the singularity is by assumption 
regular at infinity too.)  
Suppose now that there exists a continuous family in $t$ of holomorphic gauge transformations 
$h(t)=h_i(t)$ as above, written in matrix form with respect to the cyclic trivialisation of $D(0)$. 
The first row of this matrix yields functions $g_0(t,z),\ldots , g_{m-1}(t,z)$ defined for $z$ 
in a neighborhood of $p$ with $g_0(0,z)=1,g_1(0,z)=0,\ldots ,g_{m-1}(0,z)=0$, such that 
\begin{equation}\label{gauge}
   v(t,z)=g_0(t,z)w+ \cdots + g_{m-1}(t,z)w^{(m-1)}
\end{equation}
is the cyclic vector of $D(t)$. This means that there exists a continuous family in $t$ of 
polynomials $G_1(t,z), \ldots ,G_m(t,z)$ of $z$ satisfying $G_k(0,z)=G_k(z)$, 
such that the Fuchsian equation 
\begin{equation}\label{modeq}
   v^{(m)}(t,z)=\frac{G_1(t,z)}{(z-a(t))\psi(z)}v^{(m-1)}(t,z)+\cdots 
   +\frac{G_m(t,z)}{(z-a(t))^m\psi(z)^m}v(t,z)
\end{equation}
has an apparent singularity at $a(t)$. Here, we have denoted by $v^{(k)}(t,z)$ the derivative 
of order $k$ of $v(t,z)$ with respect to $z$. By (\ref{gauge}) and the assumption of holomorphicity 
of the $g_k$ in $z$, one has 
$$
    v^{(k)}(t,z)= \sum_{j=0}^k {g}^{j}_k(t,z)w^{(k-j)}(z) 
$$
for some holomorphic functions ${g}^{j}_k$ of $z$, varying continuosly in $t$. 
Then (\ref{modeq}) is equivalent to a family in $t$ of equations 
\begin{equation}\label{Heq}
   w^{(m)}(z)=\frac{H_1(t,z)}{(z-a(t))\psi(z)}w^{(m-1)}(z)+\cdots 
   +\frac{H_m(t,z)}{(z-a(t))^m\psi(z)^m}w(z), 
\end{equation} 
where $H_k(t,z)$ are holomorphic functions with respect to $z$ and continuous in $t$. 
Hence, the original vector is cyclic for $D(t)$ for any $t$, and the induced equation is 
Fuchsian in $p$ and in $a(t)$. 
Moreover, for all $|t|<\varepsilon$ and $k \in \{ 1,\ldots,m\}$
$$
   H_k(t,p)=G_k(t,p)
$$
and 
$$
   H_k(t,a(t))=G_k(t,a(t))
$$
because the ${g}^{j}_k(t,z)$ are holomorphic in $z$. Hence, these values depend continuosly on $t$. 
Notice that since this system specialises to the original (\ref{equation}) in $t=0$, and $a(0)=p$,
it follows in particular that 
\begin{equation}\label{HG}
    H_1(0,z) = (z-p)G_1(z). 
\end{equation}
By general theory of differential equations with holomorphic coefficients (see
e.g. Section 45 of \cite{for}), 
the coefficient of index $1$ of an equation evaluated at a singularity must be a strictly 
negative integer. In our situation, this means that for all $t\neq 0$ we have 
$$
    H_1(t,a(t)) \in \Z_-. 
$$
Since $H_1$ and $a(t)$ depend continuously on $t$, so does this negative integer, and 
is therefore a constant independent of $t$. But in $t=0$, (\ref{HG}) yields 
$$
    H_1(0,a(0)) = H_1(0,p) = 0,  
$$
contradicting continuity of $H_1(t,a(t))$. Therefore, one of our assumptions has to be false: 
either the family (\ref{modeq}) is not continuous in $t$, or the gauge transformations $g(t,z)$ 
are not continuous in $t$ or not holomorphic with respect to $z$; this finishes the proof.

\section{Behaviour under the Hitchin map}\label{sec:Hitchinmap}
In this section, our aim is to describe the infinitesimal effect of the Hitchin map on the space of 
Fuchsian equations. 
Since the Hitchin map is defined for moduli spaces of Higgs bundles, we first need to identify our 
moduli space $\Mod$ to a moduli space $\Mod_{Dol}$ of Higgs bundles. This is achieved by the 
non-abelian Hodge theory of \cite{sim}. 

In order to be able to describe this correspondance, we first need to recall some basic notions. 
Let $X$ be any complex projective curve and $V$ a smooth vector bundle over $X$. 
Fix a finite set of points $S$ in $X$. 
A \emph{parabolic structure} on $V$ with parabolic points in $S$ consists of a strictly decreasing 
filtration of the fiber of $V$ at the point $s$ 
$$
    V|_s = F_0V|_s \supset F_1V|_s \supset \cdots \supset F_{m_s-1}V|_s \supset  F_{m_s}V|_s = 0
$$
for each $s \in S$, together with a strictly increasing sequence of real numbers 
$$
   0 \leq \gamma^s_1 < \cdots < \gamma^s_{m_s} <1. 
$$
This notion has an obvious generalisation for any choice of half-open interval $[a,a+1[$ 
instead of $[0,1[$, that we shall also call parabolic structure. 
The parabolic structure at $s$ is said to be \emph{trivial of weight $\gamma$} if $m_s=1$ and 
$\gamma^s_1=\gamma$. 
A \emph{parabolic logarithmic connection} over $X$ with parabolic points in $S$ is a couple 
$(E,D)$, where $E$ is a holomorphic bundle over $X$ such that the underlying smooth vector 
bundle has a parabolic structure in the points of $S$, and $D$ is a logarithmic connection 
on $E$ such that for each $s \in S$ the residue $\res(s,D)$ preserves the parabolic filtration. 
A \emph{Higgs bundle} on $X$ is a couple $(\mathcal{E},\theta)$, where $\mathcal{E}$ is a 
holomorphic vector bundle over $X$, and $\theta$ is a global holomorphic 
$\End(\mathcal{E})$-valued $1$-form, called the \emph{Higgs field}. 
Finally, a \emph{parabolic Higgs bundle} on $X$ is a couple $(\mathcal{E},\theta)$ where 
$\mathcal{E}$ is as above, and $\theta$ is a global meromorphic $\End(\mathcal{E})$-valued 
$1$-form with at most simple poles in the points of $S$ and no other poles, such that 
for each $s\in S$ the residue $\res(s,\theta)$ preserves the parabolic filtration at $s$. 

Now we are able to give the identification with a moduli space of Higgs bundles stated above. 
Take $X$ to be the Riemann sphere $\P$ and set $S=P$. 
Let us endow all the logarithmic points $P$ with the trivial parabolic filtration 
and all parabolic weights equal to 
\begin{equation}\label{drparwt}
   \beta = \frac{(n-1)(m-1)}{2(n+1)}.
\end{equation}
Then, the parabolic degree of $E$ is $0$, and the moduli space $\Mod$ naturally identifies to the 
moduli space of parabolic integrable connections with the same prescribed eigenvalues of the residues 
$\{ \mu^j_1,\ldots,\mu^j_m \}$ as before, of parabolic degree $0$, with all  parabolic weights equal 
to $\beta$ and with trivial parabolic filtration. According to the Main Theorem of \cite{sim}, 
this space is then real-analytically isomorphic to the moduli space $\Mod_{Dol}$ of 
S-equivalence classes of parabolic Higgs bundles with parabolic points at $P$, 
with fixed eigenvalues of the residue of the Higgs field at any $p_j$ equal to 
\begin{equation} \label{dolres}
   \lambda^j_k = \frac{\mu^j_k-\beta}{2} \hspace{1cm} k\in \{ 1,\ldots ,m\}
\end{equation}
and fixed parabolic weights at $p_j$ equal to 
\begin{equation} \label{dolparwt}
   \alpha^j_k =\Re(\mu^j_k)-[\Re(\mu^j_k)] \hspace{1cm} k\in \{ 1,\ldots ,m\}, 
\end{equation}
where $[.]$ stands for integer part. Denote by 
\begin{equation}\label{nah}
   \Phi :\Mod \longrightarrow \Mod_{Dol}
\end{equation}
this real-analytic isomorphism.  

Let $K$ denote the canonical sheaf on $\P$. 
Denote the characteristic polynomial of the Higgs field by 
\begin{equation}\label{charpolytheta}
    \chi_{\theta}(\lambda)=\lambda^m+a_1\lambda^{m-1}+\cdots +a_m, 
\end{equation}
where $a_k\in H^0(\P,K(P)^{k})$. The Hitchin map is defined by 
\begin{align*}
   H:\Mod & \to A\subset \oplus_{k=1}^m H^0(\P,K(P)^k)\\
     (\dbar^{\E},\theta) & \mapsto (a_1,\ldots , a_m),
\end{align*}
where $A$ is the affine space $\im(H)$, determined by the locally isomonodromic requirement. 
It is classical (see Proposition 8.3 of \cite{mar} for the non-parabolic case) 
that the dimension of the fibers of this map equals half the dimension of $\Mod$, and as a 
consequence, the same relation holds for $\dim(A)$. 
The statement in the parabolic case can be proved using a Serre duality-argument 
similar to that of Subsection \ref{subsec:construction}. 

Define the composed map $H\circ \Phi$ mapping an integrable connection coming from a Fuchsian 
equation to the coefficients of the characteristic polynomial of the Higgs bundle corresponding to it. 

We then have: 

\begin{thm}\label{thm:Hitchinmap2}
For any logarithmic connection $(E_{\L},D_{\L})\in \Mod$ induced by a Fuchsian equation (\ref{equation}) 
satisfying Condition \ref{cond}, the differential of $H\circ \Phi$ restricted to the tangent space at 
$(E_{\L},D_{\L})$ of the space $\E$ of Fuchsian equations is an isomorphism. 
\end{thm}

\begin{proof} 
First, we need to show that the tangent of the inclusion $\E \hookrightarrow \Mod$ is injective 
at any point $\L$. This is an infinitesimal version of Proposition \ref{prop:gaugeequiv}. 
In concrete terms, we need to show that if $G_k(t,z)$ are a family of polynomials of degree at most 
$k(n-1)$ for $k \in \{ 1,\ldots ,m \}$ varying smoothly with respect to $t$, and such that 
denoting by $G_k'$ the derivative of $G_k$ with respect to $t$, for at least one $k$ we have 
$G_k'(0,z)\neq 0$, then there exists no $g\in H^0(\End(E_{\L}))$ with the property that 
\begin{equation}\label{dlg}
   D_{\L}g= \begin{pmatrix}
      0 & \ldots & 0 \\
      \vdots & \ddots & \vdots \\
      G_m'(0,z)  & \ldots & G_1'(0,z)
      \end{pmatrix}.
\end{equation}
Recall that in this expression $D_{\L}$ acts on endomorphisms by the adjoint action. 
Here we only treat the case $m=2$; the generalisation to higher rank goes along the same lines as 
in the proof of Proposition \ref{prop:gaugeequiv}.  
The matrix of an endomorphism $g$ satisfying the required conditions is again lower triangular, 
with constants on the diagonal: 
$$
    g= \begin{pmatrix}
     g_{11} & 0 \\
     g_{21} & g_{22}
    \end{pmatrix}. 
$$
Condition (\ref{dlg}) is 
$$
    d^{1,0}g+\left[ \begin{pmatrix}
     g_{11} & 0 \\
     g_{21} & g_{22}
    \end{pmatrix} , 
    \begin{pmatrix}
     0 & 1 \\
     G_2(0,z) & G_1(0,z) + \psi '
    \end{pmatrix}
    \right] 
    =
    \begin{pmatrix}
     0 & 0 \\
     G_2'(0,z) & G_1'(0,z)
    \end{pmatrix}.
$$
Computing the commutator, the left-hand side is equal to
$$
   \begin{pmatrix}
     g_{21} & g_{22}-g_{11} \\
     * & *
    \end{pmatrix}.
$$
Hence, such a $g$ must satisfy $g_{22}=g_{11}$ and $g_{21}=0$. Therefore, it is a multiple of the identity. 
But then, this implies $G_1'(0,z)=G_2'(0,z)$. 

Let now $\eta$ be the harmonic $1$-form representing an element in the tangent of locally isomonodromic
deformations of the equation (\ref{equation}). By Theorem \ref{thm:main}, $\eta$ is an $L^2$ harmonic 
$(1,0)$-form. It is well-known (see \cite{hit}) that the tangent of $\Phi$ maps $\eta$ to the couple 
$$
     \left( \frac{1}{2}(\eta -\eta^*)^{0,1} ,\frac{1}{2}(\eta +\eta^*)^{1,0} \right)
$$
in Higgs coordinates. Since $\eta$ is of type $(1,0)$, this couple is equal to 
$$
     \left( -\frac{\eta^*}{2} ,\frac{\eta}{2} \right).
$$
On the other hand, the tangent map of $H$ in these coordinates is projection to the second 
component. Therefore, we have 
$$
    \d (H\circ \alpha) (\eta) = \frac{\eta}{2};
$$
whence injectivity. Surjectivity then follows by equality of dimensions. 
\end{proof}

\section{Applications}\label{sec:Applications}

\subsection{\'Etale and Lagrangian properties}

\begin{cor}\label{cor:etale}
The Hitchin map $H$ restricted to the space $\E$ of deformations of the Fuchsian equation 
(\ref{equation}) is a one-to-one cover of $A$. 
Furthermore, $\E$ is an algebraic Lagrangian submanifold of $\Mod$ with respect to the holomorphic 
symplectic form $\Omega_I=-\omega_J+i\omega_K$, where $I$ and $J$ are the de Rham and Dolbeault 
complex structures on $\Mod$, and $K=IJ$. 
\end{cor}
\begin{rk}
The Lagrangian property is proved independently by Jonathan Aidan in the recent work \cite{aidan}. 
His proof relies on matrix commutator computations, hence it is entirely different from ours. 
\end{rk}
\begin{proof}
By Theorem \ref{thm:Hitchinmap2}, the map $H|_{\E}$ is \'etale. 
On the other hand, notice that by Propositions 6.13-14 of \cite{vdps} 
(see also Theorem 1.3 of \cite{Sz}) 
$\E$ is the subspace of $\Mod$ where the type of the underlying vector bundle $E$ is equal to 
(\ref{holbdl}). 
We deduce that $\E$ is an algebraic subvariety of $\Mod$; in particular, the inclusion 
$$
    \E \hookrightarrow \Mod
$$
is proper. Finiteness of $H|_{\E}$ now follows from \cite{logmar}, where it is proven that the 
Hitchin map $H$ is generically proper. 
Finally, a finite connected unramified cover of an affine space by another one is necessarily bijective. 

The tangent space of $\Mod$ is identified with endomorphism-valued $L^2$ harmonic $1$-forms. 
For such $1$-forms $\phi_1,\phi_2$ the holomorphic symplectic form of the de Rham moduli 
space can be written as 
$$
    \Omega_I(\phi_1,\phi_2) = \int_{\C} \tr (\phi_1 \wedge \phi_2). 
$$
Clearly this quantity vanishes if both $\phi_1$ and $\phi_2$ are of type $(1,0)$. 
\end{proof}


\subsection{The number of apparent singularities} 

Another application of our method is to find the smallest number $N$ such that any logarithmic 
integrable connection $D$ correspond to a Fuchsian equation with at most $N$ apparent singularities. 
By this, we mean that the equation has the same monodromy representation as $D$. 
In \cite{oht}, M. Ohtsuki proves the inequality 
\begin{equation}\label{Nlef}
    N \leq e. 
\end{equation}

For $(a_1,\ldots , a_N)\in S^N\P$ an unordered $N$-tuple of points in $\P$, let 
$\Mod_{a_1,\ldots , a_N}$ stand for the subspace of $\Mod$ composed of logarithmic 
connections $(E,D)$ that admit a corresponding Fuchsian equation with apparent singularities 
exactly in $a_1,\ldots , a_N$. 
(In case some of the $a_l$ agree, we consider Fuchsian equations with apparent singularity in 
the corresponding point with weight equal to the multiplicity of the given point, 
and in case some of the $a_l$ coincide with one of the $p_j$, the singularity must be real with the 
prescribed exponents.) 
By definition, any $(E,D)\in\Mod$ admits at least one corresponding Fuchsian equation with at most 
$N$ apparent singularities, hence $(E,D)$ lies in at least one of the $\Mod_{a_1,\ldots , a_M}$ for 
some $M\leq N$. By the above remarks, we may take $M=N$, with the additional apparent singularities 
coinciding with some of the real singularities. 
Our aim in this subsection is to prove the following result. 

\begin{cor}\label{cor:Lagrfol}
We have $N = e$, and for all $(a_1,\ldots , a_e)\in S^e\P$ the space $\Mod_{a_1,\ldots , a_e}$ 
is a Lagrangian submanifold of $\Mod$. Furthermore, the $\Mod_{a_1,\ldots , a_e}$ for different 
$(a_1,\ldots , a_e)$ are either disjoint or equal. In particular, the $\Mod_{a_1,\ldots , a_e}$ 
define a Largangian foliation of $\Mod$ integrating the subspace distribution $H^{1,0}$ of 
$T\Mod$ defined in Section \ref{sec:Hodgestr}. 
\end{cor}

\begin{rk}
Part of this statement is already contained in the article \cite{dubma} of B. Dubrovin and 
M. Mazzocco, where it is proved that the positions of the apparent singularities are part 
of a Darboux coordinate chart on the moduli space $\Mod$. 
Their result is stronger than ours in that it gives dual coordinates as well, 
but they don't show the statement about the foliation. 
Their method of proof is different from ours. 

On the other hand, the formal dimension computation of $\E_{a_1,\ldots , a_N}$ in the proof 
below is also performed in Remark 6.33 of \cite{vdps}. 
\end{rk}
\begin{proof}
In view of (\ref{Nlef}), for the first statement it is sufficient to prove $N \geq e$.
Suppose that for an $N$-tuple of points $(a_1,\ldots , a_N)$ the space $\Mod_{a_1,\ldots , a_N}$ 
is not empty, i.e. there exists at least one logarithmic connection $(E,D)$ with a corresponding 
Fuchsian equation with the same local monodromies as (\ref{equation}) and apparent singularities 
exactly in $a_1,\ldots , a_N$. Then $\Mod_{a_1,\ldots , a_N}$ maps naturally to the affine variety 
$\E_{a_1,\ldots , a_N}$ of Fuchsian equations with the prescribed local monodromy in all singular points 
$(p_0,\ldots,p_n)$ and apparent singularities in $a_1,\ldots,a_N$. 
Indeed, up to conjugacy the residue of the integrable connection induced by any such Fuchsian equation 
around $p_j$ is by definition equal to that of $D_{\L}$, and around $a_i$ is equal to the residue of an 
integrable connection with apparent singularity in that point, 
that is one whose monodromy is the identity, and the set of logarithms of the identity is discrete. 
Conversely, by Deligne's Riemann-Hilbert correspondance, any Fuchsian equation with local monodromies 
in all real singular points $p_0,\ldots,p_n$ as prescribed, and apparent singularities exactly in the 
fixed points $a_1,\ldots , a_N$ extends to a vector bundle of arbitrary degree over $\P$, hence 
defines at least one integrable connection in $\Mod$. One obtains all such extensions of the same 
degree by various choices of logarithms of the monodromy matrices of the equation. 
Since in $\Mod$ the residues of the connection at $P$ are prescribed, in order to obtain an 
extension of an element of $\E_{a_1,\ldots , a_N}$ to an element of $\Mod$, there is no choice. 
We deduce that $\Mod_{a_1,\ldots , a_N}$ is isomorphic to $\E_{a_1,\ldots , a_N}$. 

Let us compute the dimension of this space. Set 
$$
   \psi_{a_1,\ldots ,a_N}(z)=\prod_{j=1}^n(z-p_j)\prod_{i=1}^N(z-a_i).
$$
As we have a total of $n+N$ singular points at finite distance, the degree of the coefficient 
$H_k$ of a given Fuchsian equation 
\begin{equation}\label{appeq}
   w^{(m)}=\frac{H_1(z)}{\psi_{a_1,\ldots ,a_N}}w^{(m-1)}+\cdots +\frac{H_m(z)}{\psi_{a_1,\ldots ,a_N}^m}w
\end{equation}
is by Fuchs' condition smaller than or equal to $k(n+N-1)$. This gives us a total of 
$$
   m + \frac{m(m+1)(n+N-1)}{2}
$$
parameters (the coefficients of the polynomials $H_k$). 
Fixing the local monodromy in the real singular points $P$ imposes $m$ conditions on these parameters 
at each singular point, hence a total of $(n+1)m$ conditions. 
However, one of these conditions is redundant because of the residue theorem. 
This leaves us $(n+1)m-1$ independent conditions on these parameters. 
On the other hand, by Section 16.4 of \cite{ince}, the requirement that 
$a_1,\ldots , a_N$ be apparent singularities imposes $m(m+1)/2$ conditions at each $a_i$. 
Following \cite{dubma}, we call an apparent singularity \emph{special} if the exponents of 
the equation at this point are 
$$
    m,m-2,m-3,\ldots ,1,0. 
$$
\begin{prop}
There are exactly $Nm(m+1)/2$ additional independent conditions imposed on the 
coefficients of an equation by the requirement that the singularities at $a_1,\ldots , a_N$ be 
special apparent. Furthermore, these additional conditions are all independent from the previously 
obtained $(n+1)m-1$ conditions imposed by the requirement that the exponents at $p_0,\ldots ,p_n$ be fixed. 
\end{prop}
\begin{rk}
The proof can presumably be modified to a proof of the general (not necessarily special) case as well, 
but the formulae become more complicated. 
\end{rk}
\begin{proof}
Let us write 
\begin{equation}\label{H}
    H_k(z)=H_k^0+H_k^1z+\cdots + H_k^{k(N+n-1)}z^{k(N+n-1)}, 
\end{equation}
for constants $H_k^l$. 
Following \cite{ince}, let us spell out explicitly the conditions on the $H_k^l$. 

First, let us treat the case of a real singularity $p_j$. 
Here, the only conditions relate to the exponents, that is to say the roots of the polynomial in $\rho$ 
\begin{align}
   [\rho]_m  - \frac{H_1(p_j)}{\prod_{\alpha \neq j}(p_j-p_{\alpha})\prod_{i}(p_j-a_i)}[\rho]_{m-1} - \cdots 
             - \frac{H_{m}(p_j)}{\prod_{\alpha \neq j}(p_j-p_{\alpha})^m\prod_{i}(p_j-a_i)^m}, 
\end{align}
where we have written
$$
    [\rho]_l=\rho(\rho-1)\cdots (\rho-l+1). 
$$
Fixing the solutions of this equation amounts to fixing all the coefficients 
\begin{align}\label{indeqconstraints:p_j}
     H_1(p_j) & = H_1^0+H_1^1p_j+\cdots + H_1^{N+n-1}p_j^{N+n-1} \notag \\ 
              & \vdots \\
     H_m(p_j) & = H_m^0+H_m^1p_j+\cdots + H_m^{m(N+n-1)}p_j^{m(N+n-1)}. \notag
\end{align}
As these equations involve pairwise distinct sets of unknowns $H_{\alpha}^{\beta}$, they are clearly 
linearly independent among each other for a fixed $p_j$. 
Let us now consider the singularity at infinity $p_0$. 
Let $v=z^{-1}$ be a holomorphic coordinate centered at $p_0$. 
Taking into account that 
$$
      v \frac{\partial}{\partial v} = -z \frac{\partial}{\partial z},
$$
the differential equation can be rewritten as 
$$
    v^m \frac{\partial^m}{\partial v^m}w + 
    \frac{H_1(v^{-1})v^{-1}}{\psi_{a_1,\ldots ,a_N}(v^{-1})} v^{m-1} \frac{\partial^{m-1}}{\partial v^{m-1}}w - 
      \cdots + (-1)^{m} \frac{H_m(v^{-1})v^{-m}}{\psi_{a_1,\ldots ,a_N}(v^{-1})^m}w. 
$$
Develop the coefficients of $v^k \frac{\partial^k}{\partial v^k}w$ into series in $v$: 
\begin{align*}
     \frac{H_k(z)z^k}{\psi_{a_1,\ldots ,a_N}(z)^k} & = 
     \frac{H_k^0z^k+\cdots + H_k^{k(N+n-1)}z^{k(N+n)}}{(z-p_1)^k\cdots (z-p_n)^k(z-a_1)^k\cdots (z-a_N)^k} \\
     & = H_k^{k(N+n-1)} + O(v). 
\end{align*}
Hence, the indicial equation at infinity is 
\begin{equation}\label{indeqinf}
     [\rho]_m  + H_1^{N+n-1}[\rho]_{m-1} - \cdots + (-1)^{m} H_m^{m(N+n-1)}. 
\end{equation}
Therefore, fixing all the roots of this polynomial amounts to fixing all coefficients: 
\begin{align}\label{indeqconstraints:inf}
     H_1^{N+n-1} & = * \notag  \\ 
     \vdots & \\
     H_m^{m(N+n-1)} & = * \notag 
\end{align}
Let us now relate the constraints (\ref{indeqconstraints:p_j}-\ref{indeqconstraints:inf}) with each 
other for all real singularities $p_j$. 
Considering these equations for all $p_j$ and rearranging the ones 
involving the same unknowns, one is led to the set of systems: 
\begin{align}\label{indeqconstraints:p}
    H_1^0+H_1^1p_1+\cdots + H_1^{N+n-1}p_1^{N+n-1} & = * \notag \\ 
                                                 & \vdots \notag \\
    H_1^0+H_1^1p_n+\cdots + H_1^{N+n-1}p_n^{N+n-1} & = * \notag \\ 
    H_1^{N+n-1} & = * \notag  \\ 
                                                 & \vdots \\
                                                 & \vdots \notag \\
    H_m^0+H_m^1p_1+\cdots + H_m^{m(N+n-1)}p_1^{m(N+n-1)} & = * \notag \\
                                                 & \vdots \notag \\
    H_m^0+H_m^1p_n+\cdots + H_m^{m(N+n-1)}p_n^{m(N+n-1)} & = * \notag \\
    H_m^{m(N+n-1)} & = * \notag  
\end{align}
where the $*$ on the right-hand side are some explicit constants determined by the fixed exponents. 
This decomposes into $m$ systems, each one consisting of $n$ equations of the van der Monde type, 
and a further relation equating the top coefficient $H_k^{k(N+n-1)}$ of $H_k$ to a constant. 

Let us now treat the case of apparent singularities.  
At each such point, the exponents are supposed to be the fixed integers $m,m-2,\ldots ,1,0$. 
Therefore, similarly to the case of a real singularity, we obtain the conditions 
\begin{align}\label{indeqconstraints:a}
    H_1^0+H_1^1a_1+\cdots + H_1^{N+n-1}a_1^{N+n-1} & = * \notag \\ 
                                                 & \vdots \notag \\
    H_1^0+H_1^1a_N+\cdots + H_1^{N+n-1}a_N^{N+n-1} & = * \notag \\ 
                                                 & \vdots \\
                                                 & \vdots \notag \\
    H_m^0+H_m^1a_1+\cdots + H_m^{m(N+n-1)}a_1^{m(N+n-1)} & = * \notag \\
                                                 & \vdots \notag \\
    H_m^0+H_m^1a_N+\cdots + H_m^{m(N+n-1)}a_N^{m(N+n-1)} & = * \notag, 
\end{align}
where the $*$ on the right-hand side are some explicit constants determined by the exponents. 
In order to write the remaining conditions, fix $i \in \{ 1,\ldots ,N\}$ and expand (\ref{appeq}) 
around $a_i$ as 
\begin{equation}\label{appeq:a_i}
   w^{(m)}=\frac{H_{1,i}(z)}{(z-a_i)}w^{(m-1)}+\cdots +\frac{H_{m,i}(z)}{(z-a_i)^m}w, 
\end{equation}
where the $H_{k,i}$ are analytic functions of $z$: 
\begin{align}\label{H_ki}
    H_{k,i}(z)& =(z-a_i)^k\frac{H_k(z)}{\psi_{a_1,\ldots ,a_N}^k} \notag \\ 
             & =\sum_{l=0}^{\infty}H_{k,i}^l(z-a)^l. 
\end{align}
Here the coefficients are determined as 
\begin{equation}\label{H_kil}
   H_{k,i}^l = \frac{1}{l!}\frac{\d^l (H_k(z)\prod_j(z-p_j)^{-k}\prod_{\beta\neq i}(z-a_{\beta})^{-k})}{\d z^l}(a_i).
\end{equation}
For $l\geq 1$ set now 
\begin{equation}\label{f_il}
   f_{l,i}(\rho)= H_{1,i}^l[\rho]_{m-1} + \cdots + H_{m-1,i}^l \rho + H_{m,i}^l, 
\end{equation}
and form the matrices  
\begin{align}
    F_{\nu,i}(\rho) = \begin{pmatrix}
        f_{1,i}(\rho + \nu -1) & f_{2,i}(\rho + \nu -2) & \ldots & f_{\nu-1,i}(\rho +1) & f_{\nu,i}(\rho)\\
        f_{0,i}(\rho + \nu -1) & f_{1,i}(\rho + \nu -2) & \ldots & f_{\nu-2,i}(\rho +1) & f_{\nu-1,i}(\rho)\\
        0 & f_{0,i}(\rho + \nu -2) & \ldots & f_{\nu-3,i}(\rho +1) & f_{\nu-2,i}(\rho)\\ 
        \vdots &  \vdots & \ddots & \vdots & \vdots \\
        0 & 0 & \ldots & f_{0,i}(\rho +1) & f_{1,i}(\rho)
    \end{pmatrix}. 
\end{align}
Finally, denote by 
$$
   \rho_{1,i}>\rho_{2,i}>\cdots >\rho_{m,i}
$$
the exponents at $a_i$, i.e. the roots of the polynomial $f_{0,i}$ of $\rho$. 
By Section 16.33 of \cite{ince}, for all $\mu$ from $2$ to $m$, there are additional conditions on 
the coefficients. For a fixed $\mu$, these remaining conditions are obtained as follows: 
first, 
\begin{equation}\label{remcond:1}
    \det F_{\rho_{\mu-1,i}-\rho_{\mu,i},i} (\rho_{\mu,i})=0; 
\end{equation}
second, 
\begin{equation}\label{remcond:2}
   \det F_{\rho_{\mu-2,i}-\rho_{\mu,i},i} (\rho_{\mu,i})=0
\end{equation}
to the second order; and so on, untill we obtain the $(\mu-1)$-th equation 
\begin{equation}\label{remcond:m}
    \det F_{\rho_{1,i}-\rho_{\mu,i},i} (\rho_{\mu,i})=0
\end{equation}
to order $(m-1)$. Let us rewrite these equations for a special apparent singularity. 
Since the exponents are independent of $a_i$ (because they are supposed to be equal to 
exponents of a special apparent singularity for all $i$), we drop the index $i$ from $\rho_{l,i}$. 
We need to write a series of equations indexed by two integers $\mu$ and $\nu$. 
\begin{enumerate} 
\item Case $\mu=2$. Here, the only sub-case is $\nu=\rho_1-\rho_2=2$, and the only condition is 
(\ref{remcond:1}), which writes explicitely 
\begin{equation}\label{remcond:m2n2}
    \det
            \begin{pmatrix} 
                f_{1,i}(m-1) & f_{2,i}(m-2) \\
                f_{0,i}(m-1) & f_{1,i}(m-2)
            \end{pmatrix} 
         = 0. 
\end{equation}
\item Case $\mu=3$. Here, $\nu$ can take two values. 
\begin{enumerate}
\item Sub-case $\nu=\rho_2-\rho_3=1$. Here, the condition writes 
\begin{equation}\label{remcond:m3n1}
      f_{1,i}(\rho_3)=f_1(m-3)=0.
\end{equation}
\item Sub-case $\nu=\rho_1-\rho_3=3$. Here, the condition is that the determinant of $F_{3,i}(\rho_3)$ 
vanishes to order $2$. Explicitly, this matrix is 
$$
     \begin{pmatrix} 
                f_{1,i}(m-1) & f_{2,i}(m-2) & f_{3,i}(m-3) \\
                f_{0,i}(m-1) & f_{1,i}(m-2) & f_{2,i}(m-3) \\
                0 & f_{0,i}(m-2) & f_{1,i}(m-3)
            \end{pmatrix} .
$$
Now, by (\ref{remcond:m3n1}) and because $m-2$ is an exponent, the last row of this matrix vanishes. 
Also, by (\ref{remcond:m2n2}) the $2$ by $2$ minor in the upper left corner vanishes. 
Since $m-2$ is a simple exponent, the vanishing to order $2$ of the determinant is then equivalent to 
\begin{equation}\label{remcond:m3n3}
    \det 
            \begin{pmatrix} 
                f_{1,i}(m-1) & f_{3,i}(m-3) \\
                f_{0,i}(m-1) & f_{2,i}(m-3)
            \end{pmatrix} 
         = 0. 
\end{equation}
\end{enumerate}
\item Case $\mu=4$. Here, $\nu$ can take three values. 
\begin{enumerate}
\item Sub-case $\nu=\rho_3-\rho_4=1$. Here, the condition writes 
\begin{equation}\label{remcond:m4n1}
      f_{1,i}(\rho_4)=f_{1,i}(m-4)=0.
\end{equation}
\item Sub-case $\nu=\rho_2-\rho_4=2$. Here, the condition is the vanishing to order $2$ of 
the determinant of 
$$
     F_{2,i}(\rho_4)= \begin{pmatrix} 
                f_{1,i}(m-3) & f_{2,i}(m-4) \\
                f_{0,i}(m-3) & f_{1,i}(m-4)
            \end{pmatrix}.
$$
As $f_{1,i}(m-3)=0$ by (\ref{remcond:m3n1}), $f_{1,i}(m-4)=0$ by (\ref{remcond:m4n1}) and 
$f_{0,i}(m-3)=0$ exactly to order $1$ because $m-3$ is a simple exponent, 
the vanishing to order $2$ of $\det(F_{2,i}(\rho_4))$ is equivalent to 
\begin{equation}\label{remcond:m4n2}
      f_{2,i}(m-4)=0.
\end{equation}
\item Sub-case $\nu=\rho_1-\rho_4=4$. The condition is the vanishing to order $3$ of the determinant of 
$$
    F_{4,i}(\rho_4)= \begin{pmatrix} 
               f_{1,i}(m-1) & f_{2,i}(m-2) & f_{3,i}(m-3) & f_{4,i}(m-4) \\
               f_{0,i}(m-1) & f_{1,i}(m-2) & f_{2,i}(m-3) & f_{3,i}(m-4) \\
               0 & f_{0,i}(m-2) & f_{1,i}(m-3) & f_{2,i}(m-4) \\ 
               0 & 0 & f_{0,i}(m-3) & f_{1,i}(m-4)
            \end{pmatrix}.
$$
Let us develop this determinant by its last two rows. According to (\ref{remcond:m3n1}), 
(\ref{remcond:m4n1}) and (\ref{remcond:m4n2}), and because $m-2$ and $m-3$ are exponents, 
all elements in the last two rows vanish at least to order $1$. 
Therefore, the vanishing to order $3$ of $\det(F_{4,i}(\rho_4))$ is equivalent to the three conditions: 
$$
    \det \begin{pmatrix} 
            f_{1,i}(m-1) & f_{2,i}(m-2) \\
               f_{0,i}(m-1) & f_{1,i}(m-2)
            \end{pmatrix} =0,
$$
$$
    \det \begin{pmatrix} 
            f_{1,i}(m-1) & f_{3,i}(m-3) \\
               f_{0,i}(m-1) & f_{2,i}(m-3)
            \end{pmatrix} =0,
$$
and 
\begin{equation}\label{remcond:m4n4}
    \det \begin{pmatrix} 
            f_{1,i}(m-1) & f_{4,i}(m-4) \\
               f_{0,i}(m-1) & f_{3,i}(m-4)
            \end{pmatrix} =0.
\end{equation}
The first two of these are exactly equal to (\ref{remcond:m2n2}) and (\ref{remcond:m3n3}) respectively. 
Therefore, the only new condition is (\ref{remcond:m4n4}). 
\end{enumerate}
\end{enumerate}
It follows by induction that the conditions for a general $\mu \in \{2,\ldots ,m\}$ are equivalent to 
the linear equations 
\begin{equation}\label{remcond:mn1}
     f_{1,i}(m-\mu)=\cdots =f_{\mu-2,i}(m-\mu)=0
\end{equation}
and the quadratic equations 
\begin{equation}\label{remcond:mn2}
       \det \begin{pmatrix} 
            f_{1,i}(m-1) & f_{\mu,i}(m-\mu ) \\
               f_{0,i}(m-1) & f_{\mu-1,i}(m-\mu )
            \end{pmatrix} =0, 
\end{equation}
where it is understood that for $\mu=2$ the condition (\ref{remcond:mn1}) is empty. 
This makes up $\mu-1$ conditions for a fixed $\mu$, hence a total of $m(m-1)/2$ conditions. 
This is the number of the remaining equations at an apparent singularity. 

Let us now prove that the conditions obtained in (\ref{indeqconstraints:p}), (\ref{indeqconstraints:a}) 
and (\ref{remcond:mn1})-(\ref{remcond:mn2}) are independent, except for one linear relation between them. 
Let us start by showing that the sets of equations in (\ref{indeqconstraints:p}) and 
(\ref{indeqconstraints:a}) are independent of each other. 
Rearranging the equations relating to the same coefficients, one obtains the systems 
\begin{align}
    H_k^0+H_k^1p_1+\cdots + H_k^{k(N+n-1)}p_1^{k(N+n-1)} & = * \notag \\ 
                                                 & \vdots \notag \\
    H_k^0+H_k^1p_n+\cdots + H_k^{k(N+n-1)}p_n^{k(N+n-1)} & = * \notag \\ 
    H_k^0+H_k^1a_1+\cdots + H_k^{k(N+n-1)}a_1^{k(N+n-1)} & = * \notag \\ 
                                                 & \vdots \notag \\
    H_k^0+H_k^1a_N+\cdots + H_k^{k(N+n-1)}a_N^{k(N+n-1)} & = * \notag \\
    H_k^{k(N+n-1)} & = * \notag  
\end{align}
for each $1\leq k \leq m$. For $k=1$, this is an $(N+n)\times (N+n)$ system of van der Monde type
in in the unknowns $(H_1^0,\ldots H_1^{N+n-1})$, plus a further linear relation only involving $H_1^{N+n-1}$. 
As we supposed the $p_j$'s and the $a_i$'s to be different, the corresponding determinant of the 
van der Monde system is non-zero. However, the extra equation on $H_1^{N+n-1}$ is a linear combination 
of the equations of the van der Monde system, because this latter is a maximal linearly independent system. 
This means that there exists a unique relation between the constraints for $k=1$, 
and that if the constants on the right-hand side are determined so that they satisfy this linear equation, 
then the coefficients $H_1^0,\ldots ,H_1^{N+n-1}$ are uniquely determined (otherwise the system has no solution). 
In concrete terms, this relation expresses the fact that the sum of the eigenvalues of all the residues 
of the induced logarithmic connection is an integer, which is just the residue theorem. 
On the other hand, for $k\geq 2$, this system is an $(N+n)\times (k(N+n-1)+1)$ system of van der Monde type, 
plus an equation only involving $H_k^{k(N+n-1)}$. Developping the corresponding determinant according to 
its last row, we see that the independence of this system boils down to the independence of an 
$(N+n)\times k(N+n-1)$ system of van der Monde type. 
Clearly $k(N+n-1)\geq N+n$ for any $k\geq 2$ provided $n\geq 2$, so the equations are independent. 

Let us now come to the proof of the fact that the linear conditions imposed on the $H_k^l$ by 
\begin{equation}\label{lineql1}
    f_{1,i}(m-3) = f_{1,i}(m-4) = \cdots = f_{1,i}(0) =0 
\end{equation}
are independent from each other for $1\leq i \leq N$ and from the previous conditions 
(\ref{indeqconstraints:p}, \ref{indeqconstraints:a}). From equation (\ref{f_il}) for $l=1$ 
we see that (\ref{lineql1}) expresses the constraint that $m-2$ roots of a polynomial of 
degree $m$ are fixed. 
Taking into account the fact that the coefficients $H_1^0,\ldots ,H_1^{N+n-1}$ of $\rho$ of the 
order $m-1$ are uniquely determined by the previous arguments, the constraints (\ref{lineql1}) 
mean that for all $i$ the vector $(H_{2,i}^1,\ldots ,H_{m,i}^1)$ lies in a certain one-dimensional 
affine subvariety $V_i^1$ of $\mathbf{A}^{m-1}$. Let us fix an arbitrary $v_i^1\in V_i^1$. 
From equation (\ref{H_kil}) for $l=1$, it follows that 
\begin{align}\label{H_ki1}
   H_{k,i}^1 & = H_k(a_i)\frac{\d (\prod_j(z-p_j)^{-k}\prod_{\beta\neq i}(z-a_{\beta})^{-k})}{\d z}|_{z=a_i} \\
            & + \frac{H_k^1+2H_k^2a_i+\cdots +k(N+n-1)H_k^{k(N+n-1)}a_i^{k(N+n-1)-1}}
              {\prod_j(a_i-p_j)^{k}\prod_{\beta\neq i}(a_i-a_{\beta})^{k}}. 
\end{align}
By (\ref{indeqconstraints:a}), the quantity $H_k(a_i)$ is equal to some constant determined by the 
exponents $m,m-2,\ldots ,1,0$. We deduce that the first term in the formula for $H_{k,i}^1$ is equal 
to a constant only depending on $(p_1,\ldots ,p_n,a_1,\ldots ,a_N)$, hence the conditions 
(\ref{lineql1}) are equivalent to saying that the vector formed by the numbers 
$$
    H_k^1+2H_k^2a_i+\cdots +k(N+n-1)H_k^{k(N+n-1)}a_i^{k(N+n-1)-1}
$$
for $2\leq k \leq m$ is equal to some $w_i^1\in \mathbf{A}^{m-1}$. The set described by 
the $w_i^1$ as $v_i^1$ runs through $V_i^1$ is a one-dimensional affine subvariety $W_i^1$ of 
$\mathbf{A}^{m-1}$, namely a translate of $V_i^1$.  
Clearly, these conditions are independent from each other for different values of $k$, for 
they involve different sets of coefficients $H_k^l$. Let us regroup again the ones relating 
to the same sets of coefficients: 
\begin{align}\label{lineqsystl1}
    H_2^1+2H_2^2a_1+\cdots + 2(N+n-1)H_2^{2(N+n-1)}a_1^{2(N+n-1)-1} & = * \notag \\ 
                                                 & \vdots \notag \\
    H_2^1+2H_2^2a_N+\cdots + 2(N+n-1)H_2^{2(N+n-1)}a_N^{2(N+n-1)-1} & = * \notag \\ 
                                                 & \vdots \\
                                                 & \vdots \notag \\
    H_m^1+2H_m^2a_1+\cdots + m(N+n-1)H_m^{m(N+n-1)}a_1^{m(N+n-1)-1} & = * \notag \\
                                                 & \vdots \notag \\
    H_m^1+2H_m^2a_N+\cdots + m(N+n-1)H_m^{m(N+n-1)}a_N^{m(N+n-1)-1} & = * \notag, 
\end{align}
where the $*$ on the right mean some constants determined by $w_i^1$. As all the $w_i^1$ vary in $W_i^1$, 
the constants on the right-hand side form an $N$-dimensional affine variety of constants. 
This system decomposes into $m-1$ systems of size $N\times 2(N+n-1),\ldots ,N\times m(N+n-1)$ respectively, 
involving different sets of $H_k^l$'s. Furthermore, for $k$ fixed between $2$ and $m$, the set of constraints 
\begin{equation}\label{constrkl1}
    H_k^1+2H_k^2a_i+\cdots + k(N+n-1)H_k^{k(N+n-1)}a_i^{k(N+n-1)-1} = *
\end{equation}
for $1\leq i\leq N$ is obviously independent from the set of constraints 
(\ref{indeqconstraints:p}, \ref{indeqconstraints:a}), except possibly for the equations of the type 
\begin{equation}\label{constrkl0}
   H_k^0+H_k^1a_j+\cdots + H_k^{k(N+n-1)}a_j^{k(N+n-1)} = *
\end{equation}
for various $j$'s with the same $k$ or 
$$
    H_k^{k(N+n-1)} = *,
$$
because these are the only couples that involve the same sets of coefficients. 
For $k\geq 2$ fixed, the linear system formed by the equations 
(\ref{indeqconstraints:p}, \ref{indeqconstraints:a}) and (\ref{constrkl1}) for the coefficients 
$H_k^l$ has (up to a permutation of the rows) as matrix 
\begin{equation}\label{constr:indeqkl1}
     \begin{pmatrix}
         1 & p_1 & p_1^2 & \cdots & p_1^{k(N+n-1)} \\
         \vdots & \vdots & \vdots &  & \vdots \\
         1 & p_n & p_n^2 & \cdots & p_n^{k(N+n-1)} \\
         1 & a_1 &  a_1^2 &\cdots & a_1^{k(N+n-1)} \\
         \vdots & \vdots & \vdots &  & \vdots \\
         1 & a_N &  a_N^2 & \cdots & a_N^{k(N+n-1)} \\
         0 & 1 & 2a_1 &  \cdots & k(N+n-1)a_1^{k(N+n-1)-1} \\
         \vdots & \vdots & \vdots &  & \vdots \\
         0 & 1 & 2a_N &  \cdots & k(N+n-1)a_N^{k(N+n-1)-1} \\
         0 & 0 & 0 & \cdots & 1 
     \end{pmatrix}. 
\end{equation}
The size of this matrix is $(2N+n+1) \times (k(N+n-1)+1)$; hence, for all $k\geq 2$ the number of the rows is 
at least as large as the number of the columns, provided that $n\geq 2$. 
\begin{lem}\label{lem:vdM}
For any $r\geq 0$ and $b_1,\ldots ,b_{r+2}\in \C$, the determinant of 
\begin{equation}\label{bmatrix}
     \begin{pmatrix}
         1 & b_1 & b_1^2 & \cdots & b_1^{2r+1} \\
         \vdots & \vdots & \vdots &  & \vdots \\
         1 & b_{r+2} &  b_{r+2}^2 & \cdots & b_{r+2}^{2r+1} \\
         0 & 1 & 2b_1 &  \cdots & (2r+1)b_1^{2r} \\
         \vdots & \vdots & \vdots &  & \vdots \\
         0 & 1 & 2b_{r} &  \cdots & (2r+1)b_{r}^{2r} 
     \end{pmatrix} 
\end{equation}
is (up to a sign) equal to 
$$
    (b_{r+1}-b_{r+2})\prod_{1\leq i \leq r,r+1\leq j\leq r+2} (b_i-b_j)^2\prod_{1\leq i < j \leq r} (b_i-b_j)^4 . 
$$
\end{lem}
\begin{proof}
The determinant of (\ref{bmatrix}) is clearly equal to the value at $x_1=b_1,\ldots ,x_{r}=b_{r}$ of 
\begin{equation*}
     \frac{\d^{r}}{\d x_1 \cdots \d x_{r}} \det 
     \begin{pmatrix}
         1 & b_1 & b_1^2 & \cdots & b_1^{2r+1} \\
         \vdots & \vdots & \vdots &  & \vdots \\
         1 & b_{r+2} &  b_{r+2}^2 & \cdots & b_{r+2}^{2r+1} \\
         1 & x_1 & x_1^2 & \cdots & x_1^{2r+1} \\
         \vdots & \vdots & \vdots &  & \vdots \\
         1 & x_{r} &  x_{r}^2 & \cdots & x_{r}^{2r+1}
     \end{pmatrix}.
\end{equation*}
This latter is a usual van der Monde determinant, hence (\ref{bmatrix}) is equal to 
\begin{align*}
   \frac{\d^{r}\left[ \prod_{1\leq i < j \leq r+2} (b_i-b_j)\prod_{1\leq i \leq r+2,1\leq j \leq r} (b_i-x_j)
        \prod_{1\leq i < j \leq r} (x_i-x_j)\right]}{\d x_1 \cdots \d x_{r}}
\end{align*}
at $x_1=b_1,\ldots ,x_{r}=b_{r}$. Differentiating this product, the only 
term that does not vanish when setting $x_i=b_i$ is the one where each expression $(b_i-x_i)$ 
in the middle product is differentiated; that is, the term 
\begin{align*}
    \pm \prod_{1\leq i < j \leq r+2} (b_i-b_j)\prod (b_i-b_j)
        \prod_{1\leq i < j \leq r} (b_i-b_j),
\end{align*}
where the product in the middle runs over all $i \in \{ 1,\ldots ,r+2 \},i \in \{ 1,\ldots ,r \}$ 
such that $i \neq j$. Counting the powers of the various factors $(b_i-b_j)$ occuring concludes the proof. 
\end{proof}
We are now able to finish the proof of independence of the linear constraints (\ref{lineql1}) from the 
previously obtained constraints (\ref{indeqconstraints:p}, \ref{indeqconstraints:a}). Indeed, 
consider the matrix (\ref{constr:indeqkl1}) for $k=2$.  Developping its determinant 
according to the last row, we see that independence of the system of constraints is equivalent 
to non-vanishing of the determinant 
\begin{equation*}
     \begin{pmatrix}
         1 & p_1 & p_1^2 & \cdots & p_1^{2(N+n-1)-1} \\
         \vdots & \vdots & \vdots &  & \vdots \\
         1 & p_n & p_n^2 & \cdots & p_n^{2(N+n-1)-1} \\
         1 & a_1 &  a_1^2 &\cdots & a_1^{2(N+n-1)-1} \\
         \vdots & \vdots & \vdots &  & \vdots \\
         1 & a_N &  a_N^2 & \cdots & a_N^{2(N+n-1)-1} \\
         0 & 1 & 2a_1 &  \cdots & [2(N+n-1)-1]a_1^{2(N+n-1)-2} \\
         \vdots & \vdots & \vdots &  & \vdots \\
         0 & 1 & 2a_N &  \cdots & [2(N+n-1)-1]a_N^{2(N+n-1)-2} \\
     \end{pmatrix}. 
\end{equation*}
This matrix can be extended by the $n-2$ rows 
$$
     0 \quad 1 \quad 2p_j \quad  \cdots \quad [2(N+n-1)-1]p_j^{2(N+n-1)-2} 
$$ 
for $1\leq j \leq n-2$ so that the extended matrix becomes (up to a permutation of the rows) of the form 
(\ref{bmatrix}) with $r=N+n-2$ and $b_1=a_1,\ldots ,b_N=a_N,b_{N+1}=p_1,\ldots ,b_{r+2}=p_n$. 
Since these points are all supposed to be different, the lemma implies that the determinant does not vanish; 
hence the rows are linearly independent. 
For $k>2$, the same proof shows that the system of rows obtained by only keeping the first $2(N+n)-1$ 
columns is linearly independent. 

The proof of the fact that the linear conditions imposed on the $H_k^l$ by 
\begin{equation}\label{lineql2}
    f_{2,i}(m-4) = f_{2,i}(m-5) = \cdots = f_{2,i}(0) =0 
\end{equation}
are independent from each other for $1\leq i \leq N$ and from the previous conditions 
(\ref{indeqconstraints:p}, \ref{indeqconstraints:a}, \ref{lineqsystl1}) is analoguos to the 
case of (\ref{lineqsystl1}). Let us briefly sketch it. 
For all $i$, these equations fix $m-3$ roots of the polynomial 
$f_{2,i}$ of degree $m$ of the variable $\rho$. Taking into account the fact that the coefficients 
$H_1^0,\ldots ,H_1^{N+n-1}$ are uniquely determined, the sum of the remaining $3$ roots at each 
$a_i$ has to be some constant. 
Once $v_i^1\in V_i^1$ is fixed, this conditions leaves two free parameters on the roots at all $i$. 
Said differently, there exists a $2$-dimensional affine variety $V_i^2\subset \mathbf{A}^{m-2}$ 
for all $1\leq i \leq N$ (depending on the choice of $v_i^1$) such that the conditions 
(\ref{lineql2}) are equivalent to $(H_{3,i}^2,\ldots ,H_{m,i}^2)\in V_i^2$. 
Again, we fix a point $v_i^2\in V_i^2$ arbitrarily, and rewrite (\ref{H_kil}) for $k\geq 3$ as 
\begin{align*}
 H_{k,i}^2 & = H_k(a_i)\frac{\d^2 (\prod_j(z-p_j)^{-k}\prod_{\beta\neq i}(z-a_{\beta})^{-k})}{\d z^2}|_{z=a_i} \\
          & + (H_k^1+2H_k^2a_i+\cdots +k(N+n-1)H_k^{k(N+n-1)}a_i^{k(N+n-1)-1}) \\
          & \qquad \frac{\d \prod_j(z-p_j)^{-k}\prod_{\beta\neq i}(a_i-a_{\beta})^{-k}}{\d z}|_{z=a_i} \\
          & + \frac{H_k^2+6H_k^3a_i+\cdots +k(N+n-1)[k(N+n-1)-1]H_k^{k(N+n-1)}a_i^{k(N+n-1)-2}}
                   {\prod_j(a_i-p_j)^{k}\prod_{\beta\neq i}(a_i-a_{\beta})^{k}}. 
\end{align*}
By (\ref{H_ki1}), the sum of the first two terms is equal to a constant only depending on 
$p_1,\ldots ,p_n,a_1,\ldots ,a_N$ and the choice of $v_i^1\in V_i^1$. 
Hence, the new condition on the $H_{k,i}^2$ for $k\geq 3$ is that the quantity 
$$
    H_k^2+6H_k^3a_i+\cdots +k(N+n-1)[k(N+n-1)-1]H_k^{k(N+n-1)}a_i^{k(N+n-1)-2}
$$
be equal to some constant $w_i^2$ only depending on the positions of 
$(p_1,\ldots ,p_n,a_1,\ldots ,a_N)$ and the choices $v_i^1\in V_i^1,v_i^2\in V_i^2$. 
Therefore, for $k=3$ we obtain the system 
\begin{align}\label{lineq:k3}
    H_3^0+H_3^1p_1+\cdots + H_3^{3(N+n-1)}p_1^{3(N+n-1)} & = * \notag \\
                                                & \vdots \notag \\
    H_3^0+H_3^1p_n+\cdots + H_3^{3(N+n-1)}p_n^{3(N+n-1)} & = * \notag \\
    H_3^0+H_3^1a_1+\cdots + H_3^{3(N+n-1)}a_1^{3(N+n-1)} & = * \notag \\
                                                & \vdots  \\
    H_3^0+H_3^1a_N+\cdots + H_3^{3(N+n-1)}a_N^{3(N+n-1)} & = * \notag \\
    H_3^1+2H_3^2a_1+\cdots + 3(N+n-1)H_3^{3(N+n-1)}a_1^{3(N+n-1)-1} & = * \notag \\ 
                                                 & \vdots \notag \\
    H_3^1+2H_3^2a_N+\cdots + 3(N+n-1)H_3^{3(N+n-1)}a_N^{3(N+n-1)-1} & = * \notag \\ 
    H_3^2+6H_3^3a_1+\cdots +3(N+n-1)[3(N+n-1)-1]H_3^{3(N+n-1)}a_1^{3(N+n-1)-2} & = * \notag\\
                                                & \vdots  \notag \\
    H_3^2+6H_3^3a_N+\cdots +3(N+n-1)[3(N+n-1)-1]H_3^{3(N+n-1)}a_N^{3(N+n-1)-2} & = * \notag\\
    H_3^{3(N+n-1)} & = *, \notag 
\end{align}
where the $*$ on the right-hand side are some explicit constants depending on 
$p_1,\ldots ,p_n,a_1,\ldots ,a_N$ and the choices $v_i^1\in V_i^1,v_i^2\in V_i^2$. 
The matrix of this system is of size $(3N+n+1)\times [3(N+n-1)+1]$, so the number of 
rows is at least as large as the number of columns provided $n\geq 2$. 
Hence, we can extend it by $n-1$ lines of the form 
$$
     0 \quad 1 \quad 2p_j \quad  \cdots \quad [3(N+n-1)-1]p_j^{3(N+n-1)-2} 
$$ 
for $1\leq j \leq n-1$ and $n-2$ lines of the form 
$$
     0 \quad 0 \quad 1 \quad 3p_j \quad  \cdots \quad [3(N+n-1)-2][3(N+n-1)-1]p_j^{3(N+n-1)-3} 
$$
for $1\leq j \leq n-2$ so that setting $r=N+n-2$ and  
$b_1=a_1,\ldots ,b_N=a_N,b_{N+1}=p_1,\ldots ,b_{r+2}=p_n$ 
we obtain (up to a permutation of the rows) the matrix 
\begin{equation*}
     \begin{pmatrix}
         1 & b_1 & b_1^2 & \cdots & b_1^{3r+2} \\
         \vdots & \vdots & \vdots &  & \vdots \\
         1 & b_{r+2} &  b_{r+2}^2 & \cdots & b_{r+2}^{3r+2} \\
         0 & 1 & 2b_1 &  \cdots & (3r+2)b_1^{3r+1} \\
         \vdots & \vdots & \vdots &  & \vdots \\
         0 & 1 & 2b_{r+1} &  \cdots & (3r+2)b_{r+1}^{3r+1} \\
         0 & 0 & 1 & \cdots & (3r+1)(3r+2)b_1^{3r} \\
         \vdots & \vdots & \vdots &  & \vdots \\
         0 & 0 & 1 & \cdots & (3r+1)(3r+2)b_{r}^{3r}
     \end{pmatrix} 
\end{equation*}
This is again a generalized van der Monde matrix, in the sense that it can be obtained 
by differentiating some rows of a usual van der Monde matrix by the variables $b_j$ 
a certain number of times. 
Similarly to Lemma \ref{lem:vdM}, the determinant of this matrix is up to a sign 
equal to a product involving positive powers of $(b_i-b_j)$, hence non-vanishing 
in the case of our choices of $b_i$'s. This handles the case $k=3$. 

Inductively, independence of the linear equations (\ref{remcond:mn1}) can be proved using 
these ideas. The key points are that for any $k\geq 2$ we have $kN+n+1\geq k(N+n-1)+1$ 
as long as $n\geq 2$, and that the determinant of a generalized van der Monde matrix is a 
product involving positive powers of $(b_i-b_j)$. Filling the technical details is left to the reader. 

There remains to prove that adjoining the constraints imposed by the quadratic equations (\ref{remcond:mn2}) 
to the previous, the obtained system of equations in the $H_k^l$ is still independent. 
Let us consider the first of these equations, namely 
$$
    \det
            \begin{pmatrix} 
                f_{1,i}(m-1) & f_{2,i}(m-2) \\
                f_{0,i}(m-1) & f_{1,i}(m-2)
            \end{pmatrix} 
         = 0. 
$$
Since we suppose that $0,1,\ldots ,m-2,m$ are the roots of the polynomial $f_{0,i}$ of degree $m$, 
we have in particular $f_{0,i}(m-1)\neq 0$. Then vanishing of the determinant is equivalent to 
\begin{equation}\label{quadcond2}
    f_{2,i}(m-2)=\frac{f_{1,i}(m-1)f_{1,i}(m-2)}{f_{0,i}(m-1)}, 
\end{equation}
hence this is just a new condition of the same type as the ones in (\ref{lineql2}), except that we 
impose $f_{2,i}(m-2)$ to be equal to some constant determined by $v_i^1$ instead of $0$. 
However, by non-vanishing of the van der Monde determinant, imposing that the value of a polynomial of 
degree $m$ at $q\leq m$ points be equal to $q$ arbitrarily assigned constants is a set of independent 
conditions on the coefficients of the polynomial. Applying this to the polynomial 
$f_{2,i}$, we deduce that the conditions imposed on the coefficients $H_{k,i}^2$ by the equations 
(\ref{lineql2}) and (\ref{quadcond2}) is that the vector $(H_{3,i}^2,\ldots ,H_{m,i}^2)\in \tilde{V}_i^2$ 
for some $1$-dimensional affine variety $\tilde{V}_i^2$. Hence, the proof of independence of the 
conditions imposed by (\ref{lineql2}) from all the previous constraints can be repeated word by word 
to yield independence of the conditions imposed by (\ref{lineql2}) and (\ref{quadcond2}) from all the 
previous constraints. 

The proof of independence of conditions (\ref{remcond:mn2}) for higher $\mu$ from all the previous 
ones is similar, and is left to the reader. 
\end{proof}
By the proposition, there is a total of $Nm(m+1)/2$ independent conditions relating to the 
apparent singularities, which are furthermore independent of the $(n+1)m-1$ independent conditions 
imposed by the real singularities. 
Therefore, the conditions that the exponents in each of the real singularities be the same 
and that $a_1,\ldots , a_N$ be apparent singularities add up to a total of 
$$
    (n+1)m - 1 + \frac{Nm(m+1)}{2}
$$
independent conditions on the $m(m+1)(n+N-1)/2$ free parameters. It follows that the 
dimension of $\E_{a_1,\ldots , a_N}$ is 
$$
   1 - m^2 + \frac{m(m-1)(n+1)}{2}, 
$$
which is exactly $e$ (see (\ref{equality}-\ref{exactvalue})). Notice in particular that 
this number is independent of $N$. 

As $\dim(\Mod)=2e$ and the subspaces $\Mod_{a_1,\ldots , a_N}$ are of dimension $e$, 
it follows that the space parametrising the leaves must also be of dimension $e$. 
As the leaves are locally parametrised by a certain subset of $S^N\P$ divided by an equivalence 
relation, we infer $N \geq e$. This finishes the proof of the first statement. 

There remains to prove that the subspaces $\Mod_{a_1,\ldots , a_e}$ integrate the distribution $H^{1,0}$. 
Indeed, the Lagrangian property then follows in exactly the same way as in Corollary \ref{cor:etale}. 
This can be achieved similarly to the proof of Subsection \ref{subsec:char}: let 
$(E,D)\in\Mod_{a_1,\ldots , a_e}$ be arbitrary, and denote by $\L_D\in \E_{a_1,\ldots , a_e}$ 
be the corresponding Fuchsian equation with apparent singularities in the set $\{ a_1,\ldots , a_e\}$.
By the equality of dimensions, it is sufficient to show that any tangent vector 
$V \in H^{1,0}\subset T_{(E,D)}\Mod$ is tangent to a deformation of $(E,D)$ in $\Mod_{a_1,\ldots , a_e}$. 
Said differently, we need to show that the locus of apparent singularities does not change 
under local holomorphic gauge transformations. 
There are two cases to distinguish: first, the case when an apparent singularity of weight 
higher than $1$ splits into several apparent singularities of lower weight; second, the 
case where the position of an apparent singularity changes without changing weight. 
The first case can be dealt with exactly as in Subsection \ref{subsec:char}: 
continuity of the integer residue of the first coefficient yields a contradiction. 
For the second case, suppose there exists a continuous family $h(t,z)$ in $t$ of holomorphic 
matrix-valued functions defined on an open set $U$ of $\P$ containing an apparent 
singularity $a$ of the equation $\L_D$, and a family of connections $D(t)$ with 
associated equations $\L(t)$ such that
\begin{itemize}
\item $h(0,z)$ is identically equal to $\mbox{Id}$, 
\item $\L(t)$ has apparent singularity in $a(t)$, where $a(t)$ is a continuous function of 
      $t$ with $a(0)=a$ but not identically equal to $a$, 
\item the local gauge transformation $h(t)$ sends $\L(0)=\L_D$ to $\L(t)$. 
\end{itemize}
In the cyclic trivialisation of $D(t)$ over $U$ defined by $\L(t)$, we have 
$$
   D(t)=\d +\frac{A(t,z)}{z-a(t)}\d z.
$$
Moreover, 
$$
   D(t)=h(t,z).D=\d + h(t,z)^{-1}A(0,z)h(t,z)\frac{\d z}{z-a} + h(t,z)^{-1} \frac{\d h(t,z)}{\d z} \d z
$$
implies 
$$
   \frac{A(t,z)}{z-a(t)}= h(t,z)^{-1}A(0,z)h(t,z)\frac{1}{z-a} + h(t,z)^{-1} \frac{\d h(t,z)}{\d z}. 
$$
It follows that for each $t \neq 0$ the matrix-valued function 
$$
    h(t,z)^{-1} \frac{\d h(t,z)}{\d z}
$$
has a pole both in $z=a$ and $z=a(t)$. Since $h(t,z)$ is supposed to be holomorphic in $z$, 
we deduce that it must have a zero in $z=a$ and in $z=a(t)$: 
$$
   h(t,z) = (z-a)(z-a(t))\tilde{h}(t,z), 
$$
where $\tilde{h}(t,z)$ is again holomorphic in $z$. Since $h(t,z)$ is continuous in $t$ and 
$a(0)=a$, this implies that 
$$
   h(0,z)=(z-a)^2\tilde{h}(0,z), 
$$
contradicting $h(0,z)=\mbox{Id}$. 
\end{proof}

Let us say that two $e$-tuples $(a_1,\ldots ,a_e),(b_1,\ldots ,b_e) \in S^e\P$ are equivalent 
if $\Mod_{a_1,\ldots ,a_e}$ and $\Mod_{b_1,\ldots ,b_e}$ coincide. 
This raises the following: 
\begin{qn}
When are two N-tuples $(a_1,\ldots ,a_e)$ and $(b_1,\ldots ,b_e)$  equivalent? 
Are the $\Mod_{a_1,\ldots ,a_e}$ also algebraic subvarieties of $\Mod$? 
\end{qn}

\bibliographystyle{alpha}
\bibliography{Fuchsian-deformations}

\end{document}